\documentclass[oneside,leqno,11pt]{article}
\usepackage{amssymb,amsmath,latexsym}

\newcommand{\diag}{\operatorname{diag}}

\def\ge{\mathfrak{e}}

\def\gg{\mathfrak{g}}
\def\gh{\mathfrak{h}}

\def\gk{\mathfrak{k}}
\def\gl{\mathfrak{l}}

\def\go{\mathfrak{o}}

\def\gq{\mathfrak{q}}

\def\gs{\mathfrak{s}}
\def\gt{\mathfrak{t}}
\def\gu{\mathfrak{u}}

\def\gz{\mathfrak{z}}

\def\Ad{{\rm Ad}}

\def\rank{{\rm rank\,}}
\def\trace{{\rm trace\,}}

\def\res{{\rm res\,}}

\def\C{\mathbb{C}}

\def\R{\mathbb{R}}

\def\Z{\mathbb{Z}}

\newtheorem{lemma}[equation]{Lemma}

\title{Branching of Representations to Symmetric Subgroups}

\date{9 September 2009}

\author{Michael G. Eastwood\,\footnote{Research supported by the
Australian Research Council.} \; and \;
Joseph A. Wolf\,\footnote{
Research partially supported by NSF Grant DMS 99-88643, by the
Australian Research Council, and by hospitality 
from the University of Adelaide.
\endgraf
\endgraf
}
}

\begin{document}
\maketitle

\begin{abstract}
Let $\gg$ be the Lie algebra of a compact Lie group and let $\theta$ be
any automorphism of $\gg$.  Let $\gk$ denote the fixed point subalgebra 
$\gg^\theta$.  In this paper we present LiE programs
that, for any finite dimensional complex representation $\pi$ of $\gg$,
give the explicit branching $\pi|_\gk$ of $\pi$ on $\gk$.  Cases of
special interest include the cases where $\theta$ has order $2$
(corresponding to compact riemannian symmetric spaces $G/K$), where
$\theta$ has order $3$ (corresponding to compact nearly--kaehler
homogeneous spaces $G/K$), where $\theta$ has order $5$ (which include
the fascinating $5$--symmetric space $E_8/A_4A_4$), and the cases where
$\gk$ is the centralizer of a toral subalgebra of $\gg$.
\end{abstract}

\section{Introduction} \label{sec1}
\setcounter{equation}{0}

There are many situations where one wants to see the explicit
branching of a particular representation from the Lie algebra 
$\gg$ of a compact Lie group to a Lie subalgebra $\gk$.  In
many cases the situation corresponds to a compact homogeneous 
space $G/K$ of some geometric interest, such as the cases where
$G/K$ is a riemannian symmetric space, a nearly--kaehler manifold,
or the compact group realization of a complex flag manifold.  
Most cases of geometric interest have the interesting property 
that $\gk$ is the fixed point set of an automorphism $\theta$
of $\gg$.  In essentially all cases one can compute the branching 
by hand, but the time and effort involved may be extreme.  This 
situation is greatly ameliorated by use of the public domain computer 
program LiE \cite{LiE}.  In this paper we produce the LiE routines that
carry out the branching of representations from $\gg$ to $\gk$ 
explicitly when $\gk$ is the fixed point set of an automorphism
$\theta$ of $\gg$.  

One might expect the built--in branch routine of LiE to do the job for
us without any additional programming.  The problem is that LiE mixes 
up the order of simple roots, making
iteration of branching very difficult and causing serious problems for
identifying the restriction in cases where there is a symmetry of the
Dynkin diagram of $\gk$.  Worse, on each summand of the 
restricted representation it renormalizes the restriction to the
center of $\gk$ in a complicated manner, and that causes even  more 
serious problems 
in geometric and analytic applications where negativity is needed and 
is controlled by restriction to the center of $\gk$.  Our LiE routines
specifically address and solve those problems.

We developed many of these LiE routines for use in our work \cite{EW}
on the range of the double fibration transform \cite[Chapter 14]{FHW}, 
where we need
explicit information on branching from the Levi component of a parabolic 
subgroup to its intersection with a maximal compact subgroup.  These LiE
programs are based on structural information on Lie algebras and automorphisms
to be found in \cite{BdS}, \cite{dS1}, \cite{dS2}, \cite{GW1} and \cite{GW2}.

We necessarily start out by describing use of the LiE program and
how its use varies with the properties of $(\gg,\theta)$.  Thus
in Section \ref{sec2} we indicate root orderings and their role
in computing LiE's ``restriction matrix''.  Then in Section \ref{sec2a}
we reduce questions of branching to the cases where $\gg$ is 
simple and $\gk$ is a maximal $\theta$--stable subalgebra of $\gg$,
where there are three essentially different situations.  The case where
$\gg$ is simple and $\rank \gk < \rank \gg$ is described in
Section \ref{sec2b}.  It relies on information from \cite{dS1}, \cite{dS2}
and \cite{GW1}.  The case where $\gg$ is simple, $\rank \gk = \rank \gg$
and $\gk$ is not semisimple, is the subject of Section \ref{sec2c}.
It relies on information from \cite{BdS}, \cite{GW1}, and the standard
structure theory of parabolic subgroups. Then the case where
$\gg$ is simple, $\rank \gk = \rank \gg$ and $\gk$ is semisimple,
is indicated in Section \ref{sec2d}.  This is the most delicate case,
and it depends on methods from \cite{BdS}, \cite{GW1} and \cite{GW2}.

In Section \ref{sec3} we list all cases where $\gg$ is simple,
$\gk$ is $\theta$--maximal and $\rank k < \rank \gg$.  For each of
them we describe how to find the restriction matrix and we give the
listing of a LiE program that computes branching from $\gg$ to $\gk$.
The programs are (\ref{branch_A_D}), (\ref{branch_A_B}), (\ref{branch_A_C}),
(\ref{branch_D_BB}), (\ref{branch_D4_G2}), (\ref{branch_D4_A2}),
(\ref{branch_E6_F4}) and (\ref{branch_E6_C4}).  In all but two of these,
$\theta^2 = 1$ so $G/K$ is a riemannian symmetric space, and in those 
two we have $\theta^3 = 1$.

In Section \ref{sec4} we discuss the LiE programs for the cases where 
$\gg$ is simple, $\rank \gk = \rank \gg$ and $\gk$ is not semisimple.
Those essentially are the cases where $\gg$ is simple and $\gk$ is the
centralizer of a toral subalgebra, where the LiE programs are described
in Section \ref{sec2c}

Section \ref{sec5} gives the LiE branching programs for the cases where
$\gg$ is simple, $\theta^2 = 1$ and $\rank \gk = \rank \gg$.  The programs
(\ref{branch_B_DB}), (\ref{branch_D_DD}) and (\ref{branch_C_CC}) apply 
when $\gg$ is classical.  There one has no surprises on the root orders,
but when $\gg$ is exceptional the LiE program scrambles the root order
going from $\gg$ to $\gk$.  In (\ref{branch_G2_A1A1}), (\ref{branch_F4_A1C3}), 
(\ref{branch_F4_B4}), (\ref{branch_E7_A7}), (\ref{branch_E8_D8}) and
(\ref{branch_E8_E7A1}) this is fairly straightforward,
as there is not much flexibility for the location of $\gk$ inside $\gg$. 
However, in applications \cite{EW} we need to keep track of the various
simple roots, and we must deal with the fact that there are three 
combinatorially distinct $A_1A_5$'s in $E_6$ and two essentially distinct
$A_1D_6$'s in $E_7$.  This results in more programs
than one might expect, specifically in (\ref{branch_E6_A1A5}),
(\ref{branch_E6_A5A1}), (\ref{branch_E6_A5A1a}), (\ref{branch_E7_A1D6})
and (\ref{branch_E7_D6A1}).

Section \ref{sec6} completes the results of Section \ref{sec5},
providing the LiE routines for the seven remaining
cases, those where $\gg$ is simple, $\theta^3 = 1$ or $\theta^5 = 1$, and
$\rank \gk = \rank \gg$.  These routines are (\ref{branch_G2_A2}),
(\ref{branch_F4_A2A2}), (\ref{branch_E6_A2A2A2}), (\ref{branch_E7_A2A5})
(\ref{branch_E7_A5A2}), (\ref{branch_E8_A8}) and (\ref{branch_E8_E6A2}).
There, as in the exceptional group cases of Section \ref{sec5}, we
label the simple roots of $\gk$ to minimize any departure from the
root ordering of $\gg$.

As indicated in Section \ref{sec2a}, this completes the analysis of
branching of finite dimensional irreducible representations from
the Lie algebra $\gg$ of a compact Lie group to the fixed point set 
$\gk$ of any automorphism $\theta$ of $\gg$.   

\section{Restriction Matrices and Branching in LiE}\label{sec2}
\setcounter{equation}{0}

All our LiE routines are given by files with names of the form
{\tt branch\_X\_Y.lie} where {\tt X} is the LiE designation of the
type of $\gg$, e.g. {\tt E6}, and {\tt Y} is the LiE designation of the
type of $\gk$, e.g. {\tt F4}.  They are called within the LiE program
by first reading in the file, ($>$ {\tt read  branch\_X\_Y.lie})
and then giving the command ($>$ {\tt  branch\_X\_Y(v)})
where {\tt v = $[v_1, \dots , v_n]$} represents the highest weight 
$\sum v_i\xi_i$ of an irreducible representation $\pi$ of $\gg$ to be 
branched on $\gk$.  Here the $\xi_i$ are the fundamental simple highest weights.
Note that this depends on the ordering of the simple roots $\psi_i$. LiE uses
(and therefore we use) Bourbaki order \cite{Bo}, given as follows on the 
Dynkin diagrams.

{\small
\setlength{\unitlength}{.75 mm}
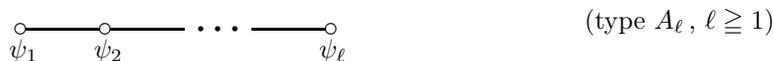
\begin{picture}(160,15)
\put(10,10){\circle{2}}
\put(8,5){$\psi_1$}
\put(11,10){\line(1,0){13}}
\put(25,10){\circle{2}}
\put(23,5){$\psi_2$}
\put(26,10){\line(1,0){13}}
\put(42,10){\circle*{1}}
\put(45,10){\circle*{1}}
\put(48,10){\circle*{1}}
\put(51,10){\line(1,0){13}}
\put(65,10){\circle{2}}
\put(63,5){$\psi_\ell$}
\put(110,10){(type $A_\ell$\,, $\ell \geqq 1$)}
\end{picture}

\setlength{\unitlength}{.75 mm}
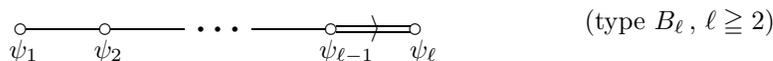
\begin{picture}(160,15)
\put(10,10){\circle{2}}
\put(8,5){$\psi_1$}
\put(11,10){\line(1,0){13}}
\put(25,10){\circle{2}}
\put(23,5){$\psi_2$}
\put(26,10){\line(1,0){13}}
\put(42,10){\circle*{1}}
\put(45,10){\circle*{1}}
\put(48,10){\circle*{1}}
\put(51,10){\line(1,0){13}}
\put(65,10){\circle{2}}
\put(63,5){$\psi_{\ell - 1}$}
\put(66,10.5){\line(1,0){13}}
\put(66,9.5){\line(1,0){13}}
\put(72,9){$\rangle$}
\put(80,10){\circle{2}}
\put(79,5){$\psi_\ell$}
\put(110,10){(type $B_\ell$\,, $\ell \geqq 2$)}
\end{picture}

\setlength{\unitlength}{.75 mm}
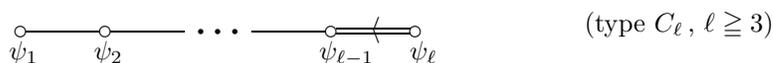
\begin{picture}(160,15)
\put(10,10){\circle{2}}
\put(8,5){$\psi_1$}
\put(11,10){\line(1,0){13}}
\put(25,10){\circle{2}}
\put(23,5){$\psi_2$}
\put(26,10){\line(1,0){13}}
\put(42,10){\circle*{1}}
\put(45,10){\circle*{1}}
\put(48,10){\circle*{1}}
\put(51,10){\line(1,0){13}}
\put(65,10){\circle{2}}
\put(63,5){$\psi_{\ell - 1}$}
\put(66,10.5){\line(1,0){13}}
\put(66,9.5){\line(1,0){13}}
\put(72,9){$\langle$}
\put(80,10){\circle{2}}
\put(79,5){$\psi_\ell$}
\put(110,10){(type $C_\ell$\,, $\ell \geqq 3$)}
\end{picture}

\setlength{\unitlength}{.75 mm}
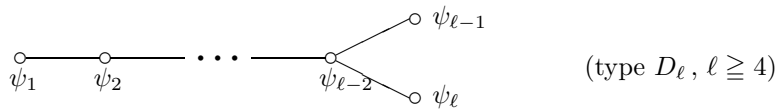
\begin{picture}(160,20)
\put(10,15){\circle{2}}
\put(8,10){$\psi_1$}
\put(11,15){\line(1,0){13}}
\put(25,15){\circle{2}}
\put(23,10){$\psi_2$}
\put(26,15){\line(1,0){13}}
\put(42,15){\circle*{1}}
\put(45,15){\circle*{1}}
\put(48,15){\circle*{1}}
\put(51,15){\line(1,0){13}}
\put(65,15){\circle{2}}
\put(63,10){$\psi_{\ell - 2}$}
\put(66,14.5){\line(2,-1){13}}
\put(80,8){\circle{2}}
\put(83,21){$\psi_{\ell-1}$}
\put(66,15.5){\line(2,1){13}}
\put(80,22){\circle{2}}
\put(83,7){$\psi_\ell$}
\put(110,12){(type $D_\ell$\,, $\ell \geqq 4$)}
\end{picture}

\setlength{\unitlength}{.75 mm}
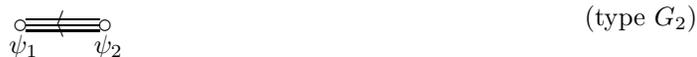
\begin{picture}(160,13)
\put(10,10){\circle{2}}
\put(8,5){$\psi_1$}
\put(11,9){\line(1,0){13}}
\put(11,10){\line(1,0){13}}
\put(11,11){\line(1,0){13}}
\put(16,9){$\langle$}
\put(25,10){\circle{2}}
\put(23,5){$\psi_2$}
\put(110,10){(type $G_2$)}
\end{picture}

\setlength{\unitlength}{.75 mm}
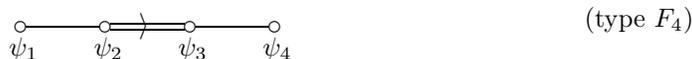
\begin{picture}(160,13)
\put(10,10){\circle{2}}
\put(8,5){$\psi_1$}
\put(11,10){\line(1,0){13}}
\put(25,10){\circle{2}}
\put(23,5){$\psi_2$}
\put(26,9.5){\line(1,0){13}}
\put(26,10.5){\line(1,0){13}}
\put(31,9){$\rangle$}
\put(40,10){\circle{2}}
\put(38,5){$\psi_3$}
\put(41,10){\line(1,0){13}}
\put(55,10){\circle{2}}
\put(53,5){$\psi_4$}
\put(110,10){(type $F_4$)}
\end{picture}

\setlength{\unitlength}{.75 mm}
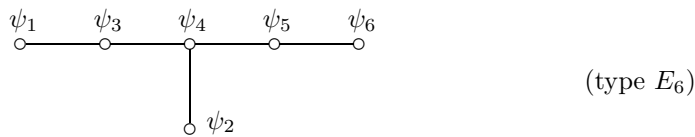
\begin{picture}(160,23)
\put(10,15){\circle{2}}
\put(8,18){$\psi_1$}
\put(11,15){\line(1,0){13}}
\put(25,15){\circle{2}}
\put(23,18){$\psi_3$}
\put(26,15){\line(1,0){13}}
\put(40,15){\circle{2}}
\put(38,18){$\psi_4$}
\put(41,15){\line(1,0){13}}
\put(55,15){\circle{2}}
\put(53,18){$\psi_5$}
\put(56,15){\line(1,0){13}}
\put(70,15){\circle{2}}
\put(68,18){$\psi_6$}
\put(40,14){\line(0,-1){13}}
\put(40,0){\circle{2}}
\put(43,0){$\psi_2$}
\put(110,7){(type $E_6$)}
\end{picture}

\setlength{\unitlength}{.75 mm}
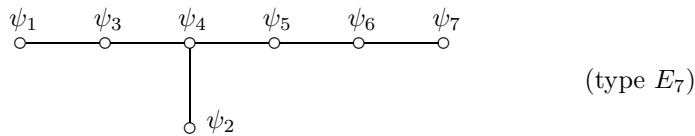
\begin{picture}(160,23)
\put(10,15){\circle{2}}
\put(8,18){$\psi_1$}
\put(11,15){\line(1,0){13}}
\put(25,15){\circle{2}}
\put(23,18){$\psi_3$}
\put(26,15){\line(1,0){13}}
\put(40,15){\circle{2}}
\put(38,18){$\psi_4$}
\put(41,15){\line(1,0){13}}
\put(55,15){\circle{2}}
\put(53,18){$\psi_5$}
\put(56,15){\line(1,0){13}}
\put(70,15){\circle{2}}
\put(68,18){$\psi_6$}
\put(71,15){\line(1,0){13}}
\put(85,15){\circle{2}}
\put(83,18){$\psi_7$}
\put(40,14){\line(0,-1){13}}
\put(40,0){\circle{2}}
\put(43,0){$\psi_2$}
\put(110,7){(type $E_7$)}
\end{picture}

\setlength{\unitlength}{.75 mm}
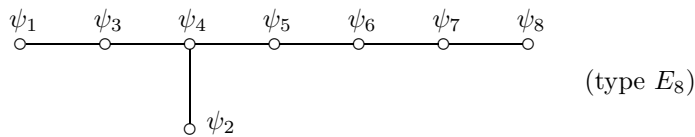
\begin{picture}(160,23)
\put(10,15){\circle{2}}
\put(8,18){$\psi_1$}
\put(11,15){\line(1,0){13}}
\put(25,15){\circle{2}}
\put(23,18){$\psi_3$}
\put(26,15){\line(1,0){13}}
\put(40,15){\circle{2}}
\put(38,18){$\psi_4$}
\put(41,15){\line(1,0){13}}
\put(55,15){\circle{2}}
\put(53,18){$\psi_5$}
\put(56,15){\line(1,0){13}}
\put(70,15){\circle{2}}
\put(68,18){$\psi_6$}
\put(71,15){\line(1,0){13}}
\put(85,15){\circle{2}}
\put(83,18){$\psi_7$}
\put(86,15){\line(1,0){13}}
\put(100,15){\circle{2}}
\put(98,18){$\psi_8$}
\put(40,14){\line(0,-1){13}}
\put(40,0){\circle{2}}
\put(43,0){$\psi_2$}
\put(110,7){(type $E_8$)}
\end{picture}
}

\noindent where, if there are two root lengths, the arrow points from
the long roots to the short roots.

If $\gk$ is a subalgebra of $\gg$ then the LiE program computes 
branching of representations by use of a ``restriction matrix''.
This is the matrix whose rows are the restrictions, from a Cartan
subalgebra $\gt$ of $\gg$ to a Cartan subalgebra $\gs \subset \gt$
of $\gk$, of the fundamental simple weights of $\gg$ as linear
combinations of the fundamental simple weights of $\gk$.  Obviously
this depends on the relation between our choices of simple root systems
for $\gg$ and $\gk$.

\subsection{Reduction to the cases where $\gg$ is simple and
$\gk$ is $\theta$--maximal.} \label{sec2a}

We start with the Lie algebra $\gg$ of a compact connected Lie group 
$G$ and an automorphism $\theta$ of $\gg$.  The fixed point algebra 
is $\gk = \gg^\theta$, and $K$ is the corresponding analytic subgroup 
of $G$.  We start also with an irreducible finite dimensional
representation $\pi$ of $\gg$.  We want to describe $\pi|_\gk$ explicitly.

We indicate how to reduce our branching questions to the case
where $\gg$ is simple and $\gk = \gg^\theta$ is maximal among the 
$\theta$--stable subalgebras of $\gg$.  That done, we have three 
essentially different possibilities.  The methods appropriate to
those three situations are addressed in Sections \ref{sec2b},
\ref{sec2c} and \ref{sec2d} below, and carried out completely in
the remainder of this paper.

Our branching procedures all use the LiE program.  We give
listings of the relevant LiE routines, and when the programming
aspects are not so obvious we give an exposition of the
mathematics behind our branching routines.

Write $\gg = \gg' \oplus \gz$ where $\gg'$ is semisimple and $\gz$ is the
center of $\gg$.  Each summand is $\theta$--stable, so
$\gk = (\gk\cap\gg')\oplus (\gk\cap\gz)$.  Also $\pi = \pi' \boxtimes \chi$,
exterior tensor product, where $\pi'$ represents $\gk\cap\gg'$ and
$\chi$ is a $1$--dimensional representation of $\gz$.  Now
$\pi|_\gk = (\pi'|_{\gk\cap\gg'})\boxtimes (\chi|_{\gk\cap\gz})$ and
evaluation of the latter factor is just restriction of a linear functional
to a linear subspace.  Thus we need only worry about computing
$\pi'|_{\gk\cap\gg'}$.  That is the first reduction: it suffices to
consider the case where $\gg$ is semisimple.

Decompose $\gg$ as a direct sum of simple ideals.  Then $\theta$ gives a
permutation on that set of ideals, and as such it is a product of
disjoint cycles.  In other words, we have a decomposition 
$\gg = \gh_1 \oplus \dots \oplus \gh_r$ where $\theta$ preserves each
$\gh_i$ and induces a cyclic permutation on its simple direct summands.
Now $\gk = \gg^\theta = \gh^\theta_1 \oplus \dots \oplus \gh^\theta_r$.
That is the second reduction: it suffices to consider
the case where $\theta$ induces a cyclic permutation on the simple ideals
of $\gg$.

Now we have reduced to the case $\gg = \gg_1 \oplus \dots \oplus \gg_m$
where the $\gg_i$ are simple, $\theta(\gg_{i-1}) = \gg_i$ for $1 < i \leqq m$,
and $\theta(\gg_m) = \gg_1$.  We interpret the $\theta: \gg_{i-1} \cong \gg_i$
as identifications.  That done, 
\begin{equation}\label{cycle1}
\gg = \gg_1 \oplus \dots \oplus \gg_1 \text{ ($m$ summands) where }
\theta(\xi_1, \dots , \xi_m) = (\gamma(\xi_m), \xi_1, \dots , \xi_{m-1})
\text{ for } \xi_i \in \gg_1 .
\end{equation}
Here $\gamma$ is an automorphism on $\gg_1$.  Now we have
\begin{equation}\label{cycle2}
\theta^m(\xi_1, \dots , \xi_m) = (\gamma(\xi_1), \dots , \gamma(\xi_m)).
\end{equation}
Thus  $\gk = \gg^\theta =
(\gg_1^\gamma \oplus \dots \oplus \gg_1^\gamma)^\theta$ where there are 
$m$ summands $\gg_1^\gamma$. Denote $\gk_1 = \gg_1^\gamma$.  From
(\ref{cycle1}) and (\ref{cycle2}) we have
\begin{equation}\label{cycle3}
\gk = \gg^\theta = 
   \{(\xi_1, \dots , \xi_1) \mid \xi_1 \in \gk_1 = \gg_1^\gamma\} 
   = \text{ diag\,}\gk_1\,.
\end{equation}
Now it suffices to consider the case where $\gg$ is simple.

We address the programming aspects.
Suppose that we are given an irreducible representation $\pi$ of
$\gg = \gg_1 \oplus \dots \oplus \gg_1$ ($m$ summands).  Then
$\pi$ is the exterior tensor product $\pi_1 \boxtimes \dots \boxtimes \pi_m$
of irreducible representation $\pi_i$ of $\gg_1$.   In view of (\ref{cycle3}),
$\pi|_\gk$ is the interior tensor product of the restrictions of the
$\pi_i$ to the $\gk_1 = \gg_1^\gamma$.  We can do this in two stages.
First we compute the restrictions $\pi_i|_{\gk_1}$, which only involves
cases where we branch from a simple Lie algebra, and then we decompose the
tensor product.  In the latter setting we have reduced to the case where
$\gamma = 1$ but $\gg_1^\gamma$ may no longer be simple.  Still, $\gg^\theta$ 
decomposes under the action of $\theta$ in the setting of a cycle of
simple ideals.  This is the third reduction: 
the branching problem is reduced to the case where $\gg = 
\gg_1 \oplus \dots \oplus \gg_1$, sum of $m$ simple ideals, and $\theta$
acts by $\theta(\xi_1, \dots , \xi_m) = (\xi_m, \xi_1, \dots , \xi_{m-1})$.
In this case $\pi = \pi_1 \otimes \dots \otimes \pi_m$ and $\gk$ is the
diagonal  diag\,$\gg_1$ in $\gg$.

We have reduced the case of branching from non--simple $\gg$ to two parts:
branching from simple proper subalgebras of $\gg$ and decomposing tensor
products of irreducible representations of $\gg$.  The latter is done in LiE
as follows.  Let {\tt v} be a matrix of $m$ rows, each row {\tt v[i]}
a vector of length equal to the rank $n$ of $\gg_1$, where the row
{\tt v[i] = [v[i,1], ... , v[i,n]]} describes the highest weight
$\lambda_i = \sum_j v[i,j]\xi_j$ of $\pi_i$ in terms of the fundamental
simple weights $\xi_j$.  If $m = 2$ and the default is set to the
Cartan type of $\gg_1$ then we can use LiE's built--in function

\centerline{\tt tensor(v[2],v[1])}

\noindent
for the tensor product decomposition.  If $m > 2$ we do this
recursively, but we must first convert the {\tt v[i]} to polynomials
in the LiE sense,

\centerline{\tt w = null(m,n); for i = i to m do w[i] = tensor(v[i],null(n)) od}

\noindent
and then we can issue the LiE command

\centerline{\tt w = w[1]; for i = 2 to m do w = tensor(w[i],w) od}

\noindent
Here is a general LiE routine to systematize this.  It is called in LiE by
{\tt branch\_diag(v,g)} where {\tt g} is the Lie type of a simple Lie
algebra such as A3, C7, G2, F4 or E8, and where {\tt v} is a matrix
of non--negative integers whose rows have length rank {\tt g} representing
highest weights of the representations of {\tt g} to be tensored together.
{\small
\begin{equation}\label{branch_diag}
\begin{array}l
\verb! # file branch_diag.lie                                          # ! \\
\verb! # usage: branch_diag(v,g) where g is a simple Lie algebra type  # ! \\
\verb! # (An, ..., E8) and v is a matrix of rank g columns, whose rows # ! \\
\verb! # specify the highest weights of reps $\pi_i$ of g;  It returns # ! \\
\verb! # the (interior) tensor product of the $\pi_i$.                 # ! \\
\verb! branch_diag(mat v; grp g) = setdefault(g);  ! \\
\verb! loc u = tensor(null(Lie_rank),null(Lie_rank)); ! \\
\verb! for r  row(v) do u = tensor(u,tensor(r,null(Lie_rank))) od; ! \\
\verb! print("the branching from product of "+n_rows(v)+" copies of " ! \\
\verb! +Lie_group(Lie_code[1],Lie_code[2])+" to the diagonal is"); u ! 
\end{array} 
\end{equation}
}

\subsection{Case $\gg$ simple and $\rank \gk < \rank \gg$.}\label{sec2b}
Suppose first that $\gg$ is simple and $\rank \gk < \rank \gg$.
Choose respective Cartan
subalgebras $\gs \subset \gt$.  Then there is a simple root system
$\Psi = \{\psi_1, \dots , \psi_n\}$ for $(\gg,\gt)$ such that the
restrictions $\psi_1|_\gs, \dots \psi_n|_\gs$ form a simple root system
$\Phi = \{\varphi_1, \dots , \varphi_r\}$ for $\gk$.  See \cite{dS1}.  In
that case we have a root restriction matrix {\tt res\_rt} whose $j^{th}$
row is given by {\tt res\_rt[j] = $[m_{j,1}, \dots m_{j,r}]$} where
$\psi_j|_\gs = \sum_k m_{j,k}\varphi_k$.  LiE however requires the
corresponding restriction matrix of fundamental simple weights, and
can compute it from {\tt res\_rt} as
\begin{equation}\label{l-rank-res_wt}
\text{\tt res\_wt = i\_Cartan($\gg$)*\res\_rt*Cartan($\gk$)/det\_Cartan($\gg$)}
\end{equation}
where {\tt i\_Cartan($\gg$)/det\_Cartan($\gg$)} is the inverse of the 
Cartan matrix of 
$\gg$ using $\Psi$ and {\tt Cartan($\gk$)} is the Cartan matrix of $\gk$
using $\Phi$.

Here is an example.  Let $\gg = \gs\gu(7)$ and $\gk = \gs\go(7)$.  The 
Dynkin diagram of $\gk$ is obtained by folding that of $\gg$,

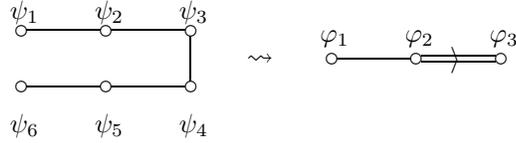
\centerline{\small
\setlength{\unitlength}{.75 mm}
\begin{picture}(100,25)
\put(10,20){\circle{2}}
\put(8,22){$\psi_1$}
\put(11,20){\line(1,0){13}}
\put(25,20){\circle{2}}
\put(23,22){$\psi_2$}
\put(26,20){\line(1,0){13}}
\put(40,20){\circle{2}}
\put(38,22){$\psi_3$}
\put(40,11){\line(0,1){8}}
\put(10,10){\circle{2}}
\put(8,2){$\psi_6$}
\put(11,10){\line(1,0){13}}
\put(25,10){\circle{2}}
\put(23,2){$\psi_5$}
\put(26,10){\line(1,0){13}}
\put(40,10){\circle{2}}
\put(38,2){$\psi_4$}
\put(50,14){$\rightsquigarrow$}
\put(65,15){\circle{2}}
\put(66,15){\line(1,0){13}}
\put(63,18){$\varphi_1$}
\put(80,15){\circle{2}}
\put(78,18){$\varphi_2$}
\put(81,15.5){\line(1,0){13}}
\put(81,14.5){\line(1,0){13}}
\put(95,15){\circle{2}}
\put(86,13.5){$\rangle$}
\put(93,18){$\varphi_3$}
\end{picture}
}
\noindent In other words, the simple root restrictions are 
$\psi_1 \mapsto \varphi_1$,
$\psi_2 \mapsto \varphi_2$,
$\psi_3 \mapsto \varphi_3$,
$\psi_4 \mapsto \varphi_3$,
$\psi_5 \mapsto \varphi_3$ and
$\psi_6 \mapsto \varphi_1$.  Thus {\tt res\_rt} is
$\left [ \begin{smallmatrix} 1 & 0 & 0 \\ 0 & 1 & 0 \\ 0 & 0 & 1 \\
 0 & 0 & 1 \\ 0 & 1 & 0 \\ 1 & 0 & 0 \end{smallmatrix} \right ]$.  Now
(\ref{l-rank-res_wt}) gives {\tt res\_wt} $= 
\left [ \begin{smallmatrix} 1 & 0 & 0 \\ 0 & 1 & 0 \\ 0 & 0 & 2 \\
 0 & 0 & 2 \\ 0 & 1 & 0 \\ 1 & 0 & 0 \end{smallmatrix} \right ]$.  If the
LiE default group is set to A6 for $\gs\gu(7)$ then branching of the adjoint
representation of $\gs\gu(7)$ on $\gs\go(7)$ is given by
{\tt branch([1,0,0,0,0,1],B3,res\_wt)}, resulting in 
{\tt 1X[0,1,0] +1X[2,0,0]}.

\subsection{Case $\gg$ simple, $\rank \gk = \rank \gg$ and $\gk$ is not 
semisimple.}\label{sec2c}

Suppose that $\gg$ is simple and $\gk$ is of equal rank but is not semisimple.
Recall that $\gk$ is $\theta$--maximal in the sense that it is maximal among
the $\theta$--stable proper subalgebras of $\gg$.  It follows that $\gk$
is the centralizer of its center, so it is a compact real form of the
reductive (Levi) component of a parabolic subalgebra of $\gg_\C$.  
We remark that the centralizer of a toral subalgebra $\ge \subset \gt$ of
$\gg$ is always the fixed point of an automorphism $\theta\in\Ad(\exp(\ge))$,
for example $\theta = \Ad(t)$ where the powers of $t$ form a dense subgroup of
the torus $\exp(\ge)$. Now $\gg$ has
a simple root system $\Psi = \{\psi_1, \dots , \psi_n\}$ such that some
subset $\Phi \subset \Psi$ is a simple root system for $\gk$.
We use the notation of Baston \& Eastwood \cite{BE}
to indicate these Levi components, i.e.\ to indicate these centralizers in
$\gg$ of subtori of its Cartan subalgebra.  Thus
if $\psi \in \Psi \setminus \Phi$ we replace the circle $\circ$  for
$\psi$ by a cross $\times$.  We refer to this as the {\em diagram} of
the corresponding parabolic subalgebra $\gq_\Phi$ of $\gg_\C$, the
corresponding parabolic subgroup $Q_\Phi$ of $G_\C$, and our algebra
$\gk = \gq_\Phi \cap \gg$ which is a compact real form of the Levi
component of $\gq_\Phi$.  For example, the parabolic
subalgebra $\gq_\Phi$ of $\gs \gl (n+1;\C)$ that corresponds to the
complex
projective space $P_n(\C) = SL(n+1;\C)/Q_\Phi$ is given by
$\Psi = \{\psi_1 , \ldots , \psi_n\}$ and
$\Phi = \{\psi_1 , \ldots , \psi_{n-1}\}$, so it has diagram
\setlength{\unitlength}{.75 mm}
\begin{picture}(52,5)
\put(3,2){\circle{2}}
\put(4,2){\line(1,0){10}}
\put(15,2){\circle{2}}
\put(16,2){\line(1,0){5}}
\put(23,2){\circle*{1}}
\put(26,2){\circle*{1}}
\put(28,2){\circle*{1}}
\put(32,2){\line(1,0){5}}
\put(38,2){\circle{2}}
\put(39,2){\line(1,0){10}}
\put(48,1){$\times$}
\end{picture},
and the parabolic subalgebra for the Grassmannian of lines
in hyperplanes in $\C^{n+1}$ is given by
$\Phi = \{\psi_2 , \ldots , \psi_{n-1}\}$ and has diagram
\setlength{\unitlength}{.75 mm}
\begin{picture}(52,5)
\put(1,1){$\times$}
\put(4,2){\line(1,0){10}}
\put(15,2){\circle{2}}
\put(16,2){\line(1,0){5}}
\put(23,2){\circle*{1}}
\put(26,2){\circle*{1}}
\put(28,2){\circle*{1}}
\put(32,2){\line(1,0){5}}
\put(38,2){\circle{2}}
\put(39,2){\line(1,0){10}}
\put(48,1){$\times$}
\end{picture}.
These correspond to the cases $\gk = \gu(n)\subset \gs\gu(n+1) = \gg$ and
$\gk = \{x \in \gu(1)\oplus\gu(n-1)\oplus\gu(1)\mid \trace x = 0\}
\subset \gs\gu(n+1) = \gg$.

Suppose that $\Phi$ consists of all but one element $\gamma = \psi_i$
of $\Psi$, in other words that $\gq_\Phi$ is a maximal parabolic subalgebra 
of $\gg_\C$.  That is the
case where there is only one $\times$ on the diagram of $\gq_\Phi$.  
Then there is a simple LiE routine (from the LiE manual~\cite{LiE}) that 
describes branching of representations
from $\gg$ to $\gk = \gq_\Phi\cap \gg$:
{\small
\begin{equation}\label{Levi_branch}
\begin{array}l
\verb! # file Levi_branch.lie # ! \\
\verb! Levi_mat(int i) = fundam(id(Lie_rank) - i) ! \\
\verb! Levi_type(int i) = Cartan_type(Levi_mat(i)) ! \\
\verb! Levi_diagram(int i) = diagram(Levi_type(i)) ! \\
\verb! Levi_res_mat(int i) = res_mat(Levi_mat(i)) ! \\
\verb! Levi_branch(vec v; int i) = loc m = Levi_mat(i); ! \\
\verb! r = res_mat(m); ! \\
\verb! branch(v, Cartan_type(m), r) ! 
\end{array}
\end{equation}}
We use it, say with $\gg = E_7$ and $\gk = E_6T_1$, as follows.  Do
{\tt read(Levi\_branch.lie)}, then
{\tt setdefault(E7)}, then {\tt diagram} in LiE to see that
$\gamma = \psi_7$, do {\tt v = [v\_1, v\_2, v\_3, v\_4, v\_5, v\_6, v\_7]}
for the highest weight $\sum v_i\xi_i$ of $\pi$, and compute the
restriction by {\tt Levi\_branch(v,7)}.
The result is a sum of vectors with multiplicities, e.g.
{\tt 1X[0,0,0,0,0,2,-6] +2X[0,0,0,0,0,2,-4] +4X[0,0,0,0,0,2,-2] + ...}
where the last entries ({\tt -6, -4, -2}) refer to the central torus. 
For the meaning of
the others do {\tt Levi\_diagram(7)} in order to compare the root
orderings (in the LiE program) for $\gg$ and $\gk$.

Suppose next that $\Phi$ consists of all but two elements $\psi_i$ and
$\psi_j$ of $\Psi$, in other words that there are two $\times$'s on the
diagram of $\gq_\Phi$.  We modify the routine (\ref{Levi_branch}) to
accommodate this.  Here it is important that $i > j$ so that we remove
rows $i$ and $j$ from a matrix by removing the $i^{\mathrm{th}}$ and then the
$j^{\mathrm{th}}$ of that.
{\small
\begin{equation}\label{Levi_branch2}
\begin{array}l
\verb! # file Levi_branch2.lie # ! \\
\verb! Levi_mat(int i, j) = fundam((id(Lie_rank) - i) - j) ! \\
\verb! Levi_type(int i, j) = Cartan_type(Levi_mat(i,j)) ! \\
\verb! Levi_diagram(int i, j) = diagram(Levi_type(i,j)) ! \\
\verb! Levi_res_mat(int i, j) = res_mat(Levi_mat(i,j)) ! \\
\verb! Levi_branch2(vec v; int i, j) = loc m = Levi_mat(i,j); ! \\
\verb! r = res_mat(m); ! \\
\verb! branch(v, Cartan_type(m), r) !
\end{array}
\end{equation}}
Similarly if $\Phi$ consists of all but three elements $\psi_i$,
$\psi_j$ and $\psi_k$ of $\Psi$, $i > j > k$, 
{\small
\begin{equation}\label{Levi_branch3.lie}
\begin{array}l
\verb! # file Levi_branch3.lie # ! \\
\verb! Levi_mat(int i, j, k) = fundam(((id(Lie_rank) - i) - j) - k) ! \\
\verb! Levi_type(int i, j, k) = Cartan_type(Levi_mat(i,j,k)) ! \\
\verb! Levi_diagram(int i, j, k) = diagram(Levi_type(i,j,k)) ! \\
\verb! Levi_res_mat(int i, j, k) = res_mat(Levi_mat(i,j,k)) ! \\
\verb! Levi_branch3(vec v; int i, j, k) = loc m = Levi_mat(i,j,k); ! \\
\verb! r = res_mat(m); ! \\
\verb! branch(v, Cartan_type(m), r)  !
\end{array}
\end{equation}}
At this point the pattern is clear.  For example, try
{\small
\begin{quote}
\begin{verbatim}
read Levi_branch3.lie
setdefault(E8)
Levi_branch3([1,0,0,0,0,0,0,1],8,6,4)
\end{verbatim}
\end{quote}}
\noindent The first five entries in each of the resulting $8$--tuples gives the 
branching on the semisimple part $A_2A_1A_1A_1$ of 
$\gk = \gq_{\Psi \setminus \{\psi_8,\psi_6,\psi_4\}}$, but with some roots
permuted.  To see the permutation look at the restriction matrix 
{\tt Levi\_res\_mat(int 8, 6, 4)} and remove rows 8, 6 and 4, and remove the
last three columns.  In LiE this can be implemented as
{\small
\begin{quote}
\begin{verbatim}
*(((*(((Levi_res_mat(8, 6, 4) - 8) - 6) - 4) -8) -7) -6)
\end{verbatim}
\end{quote}}
\noindent In this way the {\tt Levi\_branch} LiE routines give the restriction
to the semisimple part of $\gk$.  

Of course, these routines also give the action of the center of~$\gk$ on each
irreducible summand but, unfortunately, this is implemented in a rather {\em ad
hoc\/} fashion. We now explain how to specify the central action in a more
systematic and useful manner. At the same time, we avoid having to deal with
the permutations introduced by the programs \verb+Levi_branch+, as above. The
problems with these programs can be illustrated with the following simple
examples. With \verb+setdefault(F4)+ in place we have LiE calculate the
following matrices.
{\small\begin{equation}\label{matrices}\begin{tabular}{c|c|c|c}
\verb+i_Cartan+&
\verb+Levi_res_mat(3)+&\verb+Levi_res_mat(4)+&
\verb+Levi_res_mat(4,3)+\\ \hline
$\begin{array}l
\verb![[2,3,4,2]!\\
\verb!,[3,6,8,4]!\\
\verb!,[2,4,6,3]!\\
\verb!,[1,2,3,2]!\\
\verb!]!
\end{array}$&
$\begin{array}l
\verb![[0,1,0,4]!\\
\verb!,[1,0,0,8]!\\
\verb!,[0,0,0,6]!\\
\verb!,[0,0,1,3]!\\
\verb!]!
\end{array}$&
$\begin{array}l
\verb![[1,0,0,2]!\\
\verb!,[0,1,0,4]!\\
\verb!,[0,0,1,3]!\\
\verb!,[0,0,0,2]!\\
\verb!]!
\end{array}$&
$\begin{array}l
\verb![[0,1,2,0]!\\
\verb!,[1,0,4,0]!\\
\verb!,[0,0,3,0]!\\
\verb!,[0,0,0,1]!\\
\verb!]!
\end{array}$
\end{tabular}\end{equation}}
\noindent In this particular case, the matrices \verb+Levi_res_mat(i)+ are easy
to understand. The first three columns specify a permutation of the uncrossed
nodes and the last column is the $i^{\mathrm{th}}$ column of the inverse Cartan
matrix \verb+i_Cartan+. It is easy to check that the element of the Cartan
subalgebra $\gt\subset\gg$ defined by the $i^{\mathrm{th}}$ column of the
inverse Cartan matrix with respect to the basis of fundamental weights is in
the center of the corresponding Levi subalgebra~$\gk$. (Indeed, this is minus
the so-called `grading element' of the corresponding maximal parabolic
subalgebra~\cite{CS}.) Thus, the restriction matrix specifies a permutation of
the uncrossed nodes and a particular element of the center. Here is the
branching of the adjoint representation given by
\verb+Levi_branch([1,0,0,0],3)+.
$$\begin{array}r[0,1,0,4]\oplus[0,0,1,3]\oplus[1,0,2,2]\oplus[0,1,1,1]
\oplus\big([0,0,2,0]\oplus[1,1,0,0]\oplus[0,0,0,0]\big)\qquad\\[3pt]
{}\oplus[1,0,1,-1]\oplus[0,1,2,-2]\oplus[0,0,1,-3]\oplus[1,0,0,-4],
\end{array}$$
(where the ordering is given by the value of the grading element from
$-4$ to~$4$). In other words, the Lie algebra $\gg=F_4$ decomposes as
$$\gg=\gg_{-4}\oplus\gg_{-3}\oplus\gg_{-2}\oplus\gg_{-1}\oplus\gg_{0}
\oplus\gg_{1}\oplus\gg_{2}\oplus\gg_{3}\oplus\gg_{4}=
[0,1,0,4]\oplus[0,0,1,3]\oplus\cdots\oplus[1,0,0,-4]$$
(and this is exactly the realization of $\gg$ as the $|4|$-graded Lie algebra
corresponding to the parabolic subalgebra
$\gq_\Phi=\gg_0\oplus\cdots\oplus\gg_4$ as in \cite[Theorem~3.2.1]{CS}). 

The restriction matrix \verb+Levi_res_mat(4,3)+ in (\ref{matrices}) is more
difficult to understand. Certainly, we could use
{\small\begin{quote}\begin{verbatim}
[[0,1,4,2]
,[1,0,8,4]
,[0,0,6,3]
,[0,0,3,2]
]
\end{verbatim}\end{quote}}
\noindent as a more easily understandable restriction matrix. It is obtained by
using the $j^{\mathrm{th}}$ and $i^{\mathrm{th}}$ columns of \verb+i_Cartan+ to
replace the last two columns of~\verb-r-, a change that is easily implemented
in LiE by adding
{\small\begin{quote}\begin{verbatim}
for k = 1 to Lie_rank do r[k,Lie_rank-1] = i_Cartan[k,j] od;
for k = 1 to Lie_rank do r[k,Lie_rank] = i_Cartan[k,i] od;
\end{verbatim}\end{quote}}
\noindent as the penultimate two lines of \verb+Levi_branch2.lie+. In
comparison with (\ref{matrices}), the last two columns of
\verb+Levi_res_mat(4,3)+ are some linear combination of the appropriate columns
of the inverse Cartan matrix. Moreover, the case $\gg=F_4$ is deceptively
simple because its Cartan matrix has unit determinant. In general, because LiE
is restricted to integer arithmetic, it is only reasonable to use the
appropriate columns from \verb=i_Cartan=, as above. In particular, the grading
element will not be simply minus the sum of these columns but, in addition, one
must divide by \verb=det_Cartan=. Although the grading element takes on
integral values on the adjoint representation (from $-k$ to $k$ where $\gg$ is
$|k|$-graded by the parabolic subalgebra $\gq_\Phi$), for a general irreducible
representation its values will be rational with \verb+det_Cartan+ as
denominator. In any case, the raw instructions \verb+Levi_res_mat(i,j)+ and
\verb+Levi_res_mat(i,j,k)+ produce rather bizarre changes of basis from the
more natural normalization provided by the inverse Cartan matrix and even
\verb+Levi_res_mat(i)+ is better modified by 
{\small\begin{quote}\begin{verbatim}
for j = 1 to Lie_rank do r[j,Lie_rank] = i_Cartan[j,i] od;
\end{verbatim}\end{quote}}
\noindent to avoid spurious factors. 

For many purposes, however, it is better to write all weights as linear
combinations of the fundamental weights of~$(\gg,\gt)$ and, following~\cite{BE},
attach the resulting coefficients to the corresponding nodes of the Dynkin
diagram. In our example, the adjoint representation 
$\setlength{\unitlength}{1pt}
\begin{picture}(36,10)(-3,0)
\put(0,0){\circle{2}}
\put(10,0){\circle{2}}
\put(20,0){\circle{2}}
\put(30,0){\circle{2}}
\put(1,0){\line(1,0){8}}
\put(10,1){\line(1,0){10}}
\put(10,-1){\line(1,0){10}}
\put(21,0){\line(1,0){8}}
\put(15,0){\makebox(0,0){$\rangle$}}
\put(0,5){\makebox(0,0){$\scriptscriptstyle 1$}}
\put(10,5){\makebox(0,0){$\scriptscriptstyle 0$}}
\put(20,5){\makebox(0,0){$\scriptscriptstyle 0$}}
\put(30,5){\makebox(0,0){$\scriptscriptstyle 0$}}
\end{picture}$ decomposes as
\begin{equation}\label{F4decomposition}
\raisebox{-20pt}{\makebox[0pt]{$\setlength{\unitlength}{1pt}
\begin{picture}(40,10)(-5,0)
\put(0,0){\circle{2}}
\put(10,0){\circle{2}}
\put(20,0){\makebox(0,0){$\times$}}
\put(30,0){\circle{2}}
\put(1,0){\line(1,0){8}}
\put(10,1){\line(1,0){9}}
\put(10,-1){\line(1,0){9}}
\put(20,0){\line(1,0){9}}
\put(15,0){\makebox(0,0){$\rangle$}}
\put(0,5){\makebox(0,0){$\scriptscriptstyle 1$}}
\put(10,5){\makebox(0,0){$\scriptscriptstyle 0$}}
\put(20,5){\makebox(0,0){$\scriptscriptstyle 0$}}
\put(30,5){\makebox(0,0){$\scriptscriptstyle 0$}}
\end{picture}
\oplus
\begin{picture}(40,10)(-5,0)
\put(0,0){\circle{2}}
\put(10,0){\circle{2}}
\put(20,0){\makebox(0,0){$\times$}}
\put(30,0){\circle{2}}
\put(1,0){\line(1,0){8}}
\put(10,1){\line(1,0){9}}
\put(10,-1){\line(1,0){9}}
\put(20,0){\line(1,0){9}}
\put(15,0){\makebox(0,0){$\rangle$}}
\put(0,5){\makebox(0,0){$\scriptscriptstyle 0$}}
\put(10,5){\makebox(0,0){$\scriptscriptstyle 0$}}
\put(20,5){\makebox(0,0){$\scriptscriptstyle 0$}}
\put(30,5){\makebox(0,0){$\scriptscriptstyle 1$}}
\end{picture}
\oplus
\begin{picture}(40,10)(-5,0)
\put(0,0){\circle{2}}
\put(10,0){\circle{2}}
\put(20,0){\makebox(0,0){$\times$}}
\put(30,0){\circle{2}}
\put(1,0){\line(1,0){8}}
\put(10,1){\line(1,0){9}}
\put(10,-1){\line(1,0){9}}
\put(20,0){\line(1,0){9}}
\put(15,0){\makebox(0,0){$\rangle$}}
\put(0,5){\makebox(0,0){$\scriptscriptstyle 0$}}
\put(10,5){\makebox(0,0){$\scriptscriptstyle 1$}}
\put(20,5){\makebox(0,0){$\scriptscriptstyle -2$}}
\put(30,5){\makebox(0,0){$\scriptscriptstyle 2$}}
\end{picture}
\oplus
\begin{picture}(40,10)(-5,0)
\put(0,0){\circle{2}}
\put(10,0){\circle{2}}
\put(20,0){\makebox(0,0){$\times$}}
\put(30,0){\circle{2}}
\put(1,0){\line(1,0){8}}
\put(10,1){\line(1,0){9}}
\put(10,-1){\line(1,0){9}}
\put(20,0){\line(1,0){9}}
\put(15,0){\makebox(0,0){$\rangle$}}
\put(0,5){\makebox(0,0){$\scriptscriptstyle 1$}}
\put(10,5){\makebox(0,0){$\scriptscriptstyle 0$}}
\put(20,5){\makebox(0,0){$\scriptscriptstyle -1$}}
\put(30,5){\makebox(0,0){$\scriptscriptstyle 1$}}
\end{picture}
\oplus
\begin{array}c
\begin{picture}(30,10)
\put(0,0){\circle{2}}
\put(10,0){\circle{2}}
\put(20,0){\makebox(0,0){$\times$}}
\put(30,0){\circle{2}}
\put(1,0){\line(1,0){8}}
\put(10,1){\line(1,0){9}}
\put(10,-1){\line(1,0){9}}
\put(20,0){\line(1,0){9}}
\put(15,0){\makebox(0,0){$\rangle$}}
\put(0,5){\makebox(0,0){$\scriptscriptstyle 0$}}
\put(10,5){\makebox(0,0){$\scriptscriptstyle 0$}}
\put(20,5){\makebox(0,0){$\scriptscriptstyle -1$}}
\put(30,5){\makebox(0,0){$\scriptscriptstyle 2$}}
\end{picture}\\
\oplus\\
\begin{picture}(30,10)
\put(0,0){\circle{2}}
\put(10,0){\circle{2}}
\put(20,0){\makebox(0,0){$\times$}}
\put(30,0){\circle{2}}
\put(1,0){\line(1,0){8}}
\put(10,1){\line(1,0){9}}
\put(10,-1){\line(1,0){9}}
\put(20,0){\line(1,0){9}}
\put(15,0){\makebox(0,0){$\rangle$}}
\put(0,5){\makebox(0,0){$\scriptscriptstyle 1$}}
\put(10,5){\makebox(0,0){$\scriptscriptstyle 1$}}
\put(20,5){\makebox(0,0){$\scriptscriptstyle -2$}}
\put(30,5){\makebox(0,0){$\scriptscriptstyle 0$}}
\end{picture}\\
\oplus\\
\begin{picture}(30,10)
\put(0,0){\circle{2}}
\put(10,0){\circle{2}}
\put(20,0){\makebox(0,0){$\times$}}
\put(30,0){\circle{2}}
\put(1,0){\line(1,0){8}}
\put(10,1){\line(1,0){9}}
\put(10,-1){\line(1,0){9}}
\put(20,0){\line(1,0){9}}
\put(15,0){\makebox(0,0){$\rangle$}}
\put(0,5){\makebox(0,0){$\scriptscriptstyle 0$}}
\put(10,5){\makebox(0,0){$\scriptscriptstyle 0$}}
\put(20,5){\makebox(0,0){$\scriptscriptstyle 0$}}
\put(30,5){\makebox(0,0){$\scriptscriptstyle 0$}}
\end{picture}
\end{array}
\oplus
\begin{picture}(40,10)(-5,0)
\put(0,0){\circle{2}}
\put(10,0){\circle{2}}
\put(20,0){\makebox(0,0){$\times$}}
\put(30,0){\circle{2}}
\put(1,0){\line(1,0){8}}
\put(10,1){\line(1,0){9}}
\put(10,-1){\line(1,0){9}}
\put(20,0){\line(1,0){9}}
\put(15,0){\makebox(0,0){$\rangle$}}
\put(0,5){\makebox(0,0){$\scriptscriptstyle 0$}}
\put(10,5){\makebox(0,0){$\scriptscriptstyle 1$}}
\put(20,5){\makebox(0,0){$\scriptscriptstyle -2$}}
\put(30,5){\makebox(0,0){$\scriptscriptstyle 1$}}
\end{picture}
\oplus
\begin{picture}(40,10)(-5,0)
\put(0,0){\circle{2}}
\put(10,0){\circle{2}}
\put(20,0){\makebox(0,0){$\times$}}
\put(30,0){\circle{2}}
\put(1,0){\line(1,0){8}}
\put(10,1){\line(1,0){9}}
\put(10,-1){\line(1,0){9}}
\put(20,0){\line(1,0){9}}
\put(15,0){\makebox(0,0){$\rangle$}}
\put(0,5){\makebox(0,0){$\scriptscriptstyle 1$}}
\put(10,5){\makebox(0,0){$\scriptscriptstyle 0$}}
\put(20,5){\makebox(0,0){$\scriptscriptstyle -2$}}
\put(30,5){\makebox(0,0){$\scriptscriptstyle 2$}}
\end{picture}
\oplus
\begin{picture}(40,10)(-5,0)
\put(0,0){\circle{2}}
\put(10,0){\circle{2}}
\put(20,0){\makebox(0,0){$\times$}}
\put(30,0){\circle{2}}
\put(1,0){\line(1,0){8}}
\put(10,1){\line(1,0){9}}
\put(10,-1){\line(1,0){9}}
\put(20,0){\line(1,0){9}}
\put(15,0){\makebox(0,0){$\rangle$}}
\put(0,5){\makebox(0,0){$\scriptscriptstyle 0$}}
\put(10,5){\makebox(0,0){$\scriptscriptstyle 0$}}
\put(20,5){\makebox(0,0){$\scriptscriptstyle -1$}}
\put(30,5){\makebox(0,0){$\scriptscriptstyle 1$}}
\end{picture}
\oplus
\begin{picture}(38,10)(-5,0)
\put(0,0){\circle{2}}
\put(10,0){\circle{2}}
\put(20,0){\makebox(0,0){$\times$}}
\put(30,0){\circle{2}}
\put(1,0){\line(1,0){8}}
\put(10,1){\line(1,0){9}}
\put(10,-1){\line(1,0){9}}
\put(20,0){\line(1,0){9}}
\put(15,0){\makebox(0,0){$\rangle$}}
\put(0,5){\makebox(0,0){$\scriptscriptstyle 0$}}
\put(10,5){\makebox(0,0){$\scriptscriptstyle 1$}}
\put(20,5){\makebox(0,0){$\scriptscriptstyle -2$}}
\put(30,5){\makebox(0,0){$\scriptscriptstyle 0$}}
\end{picture}.$}}\end{equation}
The conversion between these two conventions is the definition of the
restriction matrix. Therefore, no matter what restriction matrix is used, to
convert back to the conventions of~\cite{BE}, one simply needs to invert the
restriction matrix and apply this inverse matrix by right multiplication to
each term obtained from \verb+branch(v, Cartan_type(m), r)+. Since LiE allows
only integer multiplication, inverting a matrix with integer entries requires
some care. In the following program, the restriction matrix \verb+r+ is
inverted by first extracting from it a permutation matrix \verb+p+, noting that
permutation matrices are orthogonal, and forming \verb+(*p)*r+. The result
necessarily has the form 
{\small$$\left[\begin{array}{ccccc}1&0&\cdots&0&*\\
0&1&\cdots&0&*\\
\vdots&\vdots&\ddots&\vdots&*\\
0&0&\cdots&1&*\\
0&0&0&0&*
\end{array}\right],$$}
which is inverted by a dint of an explicit formula.
{\small
\begin{equation}\label{Levi_branch_improved}
\begin{array}l
\verb! # file Levi_branch_improved.lie # ! \\
\verb! Levi_mat(int i) = fundam(id(Lie_rank) - i) ! \\
\verb! Levi_type(int i) = Cartan_type(Levi_mat(i)) ! \\
\verb! Levi_diagram(int i) = diagram(Levi_type(i)) ! \\
\verb! Levi_res_mat(int i) = res_mat(Levi_mat(i)) ! \\
\verb! div(pol p;int k) = loc l = 0X(null(n_vars(p))); ! \\
\verb! for i=1 to length(p) do l = l+coef(p,i)X(expon(p,i)/k) od; l ! \\
\verb! Levi_branch(vec v; int i) = loc m = Levi_mat(i); ! \\
\verb! loc r = res_mat(m); loc p = r; ! \\
\verb! for j = 1 to Lie_rank do p[j,Lie_rank] = 0 od; p[i,Lie_rank] = 1; ! \\
\verb! loc q = (*p)*r; loc det_q = q[Lie_rank,Lie_rank]; ! \\
\verb! loc qq = q; ! \\
\verb! for j = 1 to Lie_rank do qq[j,Lie_rank] = -q[j,Lie_rank] od; ! \\
\verb! qq[Lie_rank,Lie_rank] = 1; ! \\
\verb! loc s = null(Lie_rank,Lie_rank); ! \\
\verb! for j = 1 to Lie_rank do s[j,j] = det_q od; ! \\
\verb! s[Lie_rank,Lie_rank] = 1; loc i_q = qq*s; ! \\
\verb! loc b = branch(v, Cartan_type(m), r); div(b*i_q*(*p),det_q) !
\end{array}
\end{equation}}
The program is used as before but the result is expressed using the
diagrammatic conventions of~\cite{BE}. For example,
{\small\begin{quote}\begin{verbatim}
read Levi_branch_improved.lie
setdefault(F4)
Levi_branch([1,0,0,0],3)
\end{verbatim}\end{quote}}
\noindent gives
{\small\begin{quote}\begin{verbatim}
1X[0,0,-1,1] +1X[0,0,-1,2] +1X[0,0, 0,0] +1X[0,0, 0,1] +
1X[0,1,-2,0] +1X[0,1,-2,1] +1X[0,1,-2,2] +1X[1,0,-2,2] +
1X[1,0,-1,1] +1X[1,0, 0,0] +1X[1,1,-2,0]
\end{verbatim}\end{quote}}
\noindent as in~(\ref{F4decomposition}) (but devoid of the convenient ordering 
there. The ordering of (\ref{F4decomposition}) is essential when the
action of the full parabolic $\gq_\Phi$ is considered rather than just its Levi
factor. Representations of $\gq_\Phi$ are generally filtered. For example, the 
tail of (\ref{F4decomposition}),
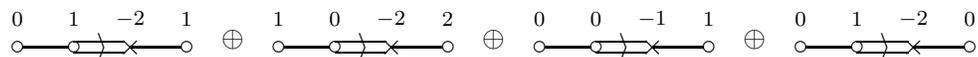
$$\begin{picture}(40,10)(-5,0)
\put(0,0){\circle{2}}
\put(10,0){\circle{2}}
\put(20,0){\makebox(0,0){$\times$}}
\put(30,0){\circle{2}}
\put(1,0){\line(1,0){8}}
\put(10,1){\line(1,0){9}}
\put(10,-1){\line(1,0){9}}
\put(20,0){\line(1,0){9}}
\put(15,0){\makebox(0,0){$\rangle$}}
\put(0,5){\makebox(0,0){$\scriptstyle 0$}}
\put(10,5){\makebox(0,0){$\scriptstyle 1$}}
\put(20,5){\makebox(0,0){$\scriptstyle -2$}}
\put(30,5){\makebox(0,0){$\scriptstyle 1$}}
\end{picture}
\oplus
\begin{picture}(40,10)(-5,0)
\put(0,0){\circle{2}}
\put(10,0){\circle{2}}
\put(20,0){\makebox(0,0){$\times$}}
\put(30,0){\circle{2}}
\put(1,0){\line(1,0){8}}
\put(10,1){\line(1,0){9}}
\put(10,-1){\line(1,0){9}}
\put(20,0){\line(1,0){9}}
\put(15,0){\makebox(0,0){$\rangle$}}
\put(0,5){\makebox(0,0){$\scriptstyle 1$}}
\put(10,5){\makebox(0,0){$\scriptstyle 0$}}
\put(20,5){\makebox(0,0){$\scriptstyle -2$}}
\put(30,5){\makebox(0,0){$\scriptstyle 2$}}
\end{picture}
\oplus
\begin{picture}(40,10)(-5,0)
\put(0,0){\circle{2}}
\put(10,0){\circle{2}}
\put(20,0){\makebox(0,0){$\times$}}
\put(30,0){\circle{2}}
\put(1,0){\line(1,0){8}}
\put(10,1){\line(1,0){9}}
\put(10,-1){\line(1,0){9}}
\put(20,0){\line(1,0){9}}
\put(15,0){\makebox(0,0){$\rangle$}}
\put(0,5){\makebox(0,0){$\scriptstyle 0$}}
\put(10,5){\makebox(0,0){$\scriptstyle 0$}}
\put(20,5){\makebox(0,0){$\scriptstyle -1$}}
\put(30,5){\makebox(0,0){$\scriptstyle 1$}}
\end{picture}
\oplus
\begin{picture}(38,10)(-5,0)
\put(0,0){\circle{2}}
\put(10,0){\circle{2}}
\put(20,0){\makebox(0,0){$\times$}}
\put(30,0){\circle{2}}
\put(1,0){\line(1,0){8}}
\put(10,1){\line(1,0){9}}
\put(10,-1){\line(1,0){9}}
\put(20,0){\line(1,0){9}}
\put(15,0){\makebox(0,0){$\rangle$}}
\put(0,5){\makebox(0,0){$\scriptstyle 0$}}
\put(10,5){\makebox(0,0){$\scriptstyle 1$}}
\put(20,5){\makebox(0,0){$\scriptstyle -2$}}
\put(30,5){\makebox(0,0){$\scriptstyle 0$}}
\end{picture}$$
is interpreted in \cite[p.~135]{BE} as inducing the cotangent bundle on the
corresponding generalized flag manifold.).

\verb+Levi_branch2.lie+ is similarly improved
{\small
\begin{equation}\label{Levi_branch2_improved}
\begin{array}l
\verb! # file Levi_branch2_improved.lie # ! \\
\verb! Levi_mat(int i, j) = fundam((id(Lie_rank) - i) - j) ! \\
\verb! Levi_type(int i, j) = Cartan_type(Levi_mat(i,j)) ! \\
\verb! Levi_diagram(int i, j) = diagram(Levi_type(i,j)) ! \\
\verb! Levi_res_mat(int i, j) = res_mat(Levi_mat(i,j)) ! \\
\verb! div(pol p;int k) = loc l = 0X(null(n_vars(p))); ! \\
\verb! for i=1 to length(p) do l = l+coef(p,i)X(expon(p,i)/k) od; l ! \\
\verb! Levi_branch2(vec v; int i, j) = loc m = Levi_mat(i,j); ! \\
\verb! loc r = res_mat(m); loc p = r; ! \\
\verb! for k = 1 to Lie_rank do p[k,Lie_rank] = 0 od; p[i,Lie_rank] = 1; ! \\
\verb! for k = 1 to Lie_rank do p[k,Lie_rank-1] = 0 od; p[j,Lie_rank-1] = 1; ! \\
\verb! loc q = (*p)*r; ! \\
\verb! loc det_q = q[Lie_rank-1,Lie_rank-1]*q[Lie_rank,Lie_rank] \ ! \\
\verb!            -q[Lie_rank,Lie_rank-1]*q[Lie_rank-1,Lie_rank]; ! \\
\verb! qq = q; ! \\
\verb! for j = 1 to Lie_rank do qq[j,Lie_rank] = -q[j,Lie_rank] od; ! \\
\verb! qq[Lie_rank-1,Lie_rank] = 0; qq[Lie_rank,Lie_rank] = 1; ! \\
\verb! for j = 1 to Lie_rank do qq[j,Lie_rank-1] = -q[j,Lie_rank-1] od; ! \\
\verb! qq[Lie_rank-1,Lie_rank-1] = 1; qq[Lie_rank,Lie_rank-1] = 0; ! \\
\verb! loc s = null(Lie_rank,Lie_rank); ! \\
\verb! for j = 1 to Lie_rank do s[j,j] = det_q od; ! \\
\verb! s[Lie_rank-1,Lie_rank-1] = (q*qq)[Lie_rank,Lie_rank]; ! \\
\verb! s[Lie_rank,Lie_rank-1] = -(q*qq)[Lie_rank,Lie_rank-1]; ! \\
\verb! s[Lie_rank-1,Lie_rank] = -(q*qq)[Lie_rank-1,Lie_rank]; ! \\
\verb! s[Lie_rank,Lie_rank] = (q*qq)[Lie_rank-1,Lie_rank-1]; ! \\
\verb! loc i_q = qq*s; ! \\
\verb! loc b = branch(v, Cartan_type(m), r); div(b*i_q*(*p),det_q) !
\end{array}
\end{equation}}
and to improve \verb+Levi_branch3.lie+ is left as an exercise (which implicitly
requires incorporating the formula for the inverse of a general
$3\times 3$ matrix).

\subsection{Case $\gg$ simple, $\rank \gk = \rank \gg$ and $\gk$ is semisimple.}
\label{sec2d}

This is the most delicate case:  $\gg$ is simple, $\theta$ is an inner
automorphism, and $\gk$ is semisimple.  Let us assume that $\gk$ is
$\theta$--maximal.  Then the Borel--de-Siebenthal structure theory 
\cite{BdS} provides a simple root system $\Psi = \{\psi_1, \dots , \psi_n\}$
and a simple root $\gamma = \psi_r \in \Psi$ 
such that $\Psi_\gk = (\Psi \setminus \{\gamma\})\cup \{-\beta_\gg\}$
is a simple root system for $\gk$, and $\theta$ has order $n_r$, where
the maximal root $\beta_\gg = \sum n_i\psi_i$.

Let $\gs$ denote the maximal rank subalgebra of $\gg$ with simple root system
$\Psi_\gs = (\Psi \setminus \{\gamma\})$.  Let $w_\gg$ and
$w_\gs$ denote, respectively, the longest elements of the Weyl groups
$W_\gg$ and $W_\gs$.  We write $\Sigma^+(\gg,\gt)$ for the positive root
system of $\gg$ relative to $\gt$ defined by $\Psi$.

\begin{lemma}\label{wg}
The transformation $-w_\gs$ preserves $\Phi_\gs$ and sends $-\beta_\gg$ 
into $\Sigma^+(\gg,\gt)$.  
\end{lemma}

\noindent {\bf Proof.}  In general, the longest element of the Weyl group
sends the positive Weyl chamber to its negative, so $-w_\gs$ preserves
$\Phi_\gs$.  But $-w_\gs(-\beta_\gg) = w_\gs(\beta_\gg)$ is obtained
from $\beta_\gg = \sum m_i \psi_i$ by a series of simple root reflections
$s_\psi : \xi \mapsto \xi - \tfrac{2\langle \psi,\xi\rangle}
{\langle \psi, \psi, \rangle}\psi$ with $\psi \ne \gamma$.  Thus the
coefficient of $\gamma$ in $w_\gs(\beta_\gg)$ is the same as that in
$\beta_\gg$, which is $n_r > 0$, so $w_\gs(\beta_\gg) \in \Sigma^+(\gg,\gt)$.  
\hfill $\square$

We now indicate how the LiE program uses $w_\gs$ and $\beta_\gg$ to compute
the restriction matrix {\tt res\_wt}, which it uses to calculate
restrictions of representations of $\gg$ to $\gk$.  
First, we use $-w_\gs$ to carry the simple root system 
$\Psi_\gk = \Psi_\gs \cup\{-\beta_\gg\}$ of $\gk$ to
another simple root system $\Phi := \Psi_\gs \cup \{w_\gs(\beta_\gg)\}$.
The point is that $\Phi$ then consists of positive roots for $\gg$, 
all but one of them simple, by
Lemma \ref{wg}.  The LiE program assumes Bourbaki root order for both $\gg$ 
and $\gk$.  It permutes the roots of $\Phi$ in a somewhat arbitrary
way in order to do this when it computes the restriction
matrix and applies it to branching of representations from $\gg$ to $\gk$.
We will try to do this in a way that involves minimal permutation.

We start by computing $w_\gs$ within the LiE program.  It is the
long word for $W_\gs$, but there LiE orders the roots incorrectly, so we
use the slightly convoluted routine
{\small
\begin{quote}
\begin{verbatim}
ws = reduce(long_word^r_reduce(long_word,[1,2,...,r-1,r+1,...,n-1,n]))
\end{verbatim}
\end{quote}}
Here {\tt long\_word} is the longest element $w_\gg$ of the Weyl group $W_\gg$, 
so the shortest element of the coset  $w_\gg W_\gs$ is
{\tt r\_reduce(long\_word,[1,2,...,r-1,r+1,...,n-1,n])}.  Then we set up 
the new simple root system $\Phi$ for $\gk$ as the rows of a matrix {\tt RR} 
in which $w_\gs(\beta_\gg)$ replaces $\gamma$ and the roots are re--ordered 
(minimally) to Bourbaki order.  Then there are three ways to compute the
restriction matrix {\tt res\_wt}.

The first is to note that {\tt RR} is the inverse of the matrix {\tt res\_rt},
so one can compute \hfill\newline
{\small\centerline{\tt res\_wt = i\_Cartan($\gg$)*\res\_rt*Cartan($\gk$)/det\_Cartan($\gg$),}}
which is (\ref{l-rank-res_wt}).  The second is just to use the LiE
assignment {\tt res\_wt = res\_mat(RR)}.  And the third, which is in fact
the way that LiE implements {\tt res\_mat}, is to set $\gg$ as the
default by {\tt setdefault($\gg$)} (putting in the Lie type of $\gg$),
initialize {\tt res\_wt}
as a the $n \times n$ identity matrix, {\tt res\_wt = id(n)}, and then fill 
it in by
{\small
\begin{quote}
\begin{verbatim}
for = 1 to n do 
   for i = 1 to n do 
      res_wt[i,j] = Cartan(i_Cartan[i],RR[j])/det_Cartan
   od 
od
\end{verbatim}
\end{quote}}
In the next few sections we will run through the various basic cases 
cases of $(\gg,\gk)$ and then reduce the general case to these basic
cases.

\section{Cases: $\gg$ is simple, $\gk$ is $\theta$--maximal and 
$\rank \gk < \rank \gg$}
\label{sec3}
\setcounter{equation}{0}
Recall the automorphism $\theta$ of $\gg$ with $\gk = \gg^\theta$.  
In this section we assume that $\gk$ is $\theta$--maximal, in
other words that it is maximal among the proper $\theta$--invariant
subalgebras of $\gg$, and we apply the methods of Section \ref{sec2b}.

Note that $\theta$ is an outer automorphism of $\gg$ because
$\rank \gk < \rank \gg$.  If some power $\theta^m \ne 1$ is an
inner automorphism then its fixed point set is $\theta$--invariant
and satisfies $\gk \subsetneqq \gg^{\theta^m} \subsetneqq \gg$.
As $\gk$ is $\theta$--maximal we conclude that every power
$\theta^m \ne 1$ is an outer automorphism of $\gg$.  All possibilities
are listed in \cite[Theorem 5.10(3)]{GW1}.  There $\theta$ has prime
order $p = 2$ or $p = 3$.  If $p = 2$ then $G/K$ is one of the 
riemannian symmetric spaces 
$$
SU(n)/SO(n),\, SU(2n)/Sp(n),\,
SO(2p+2+2q)/\{SO(2p+1)\times SO(2q+1)\},\,
E_6/F_4,\, E_6/Sp(4).
$$  
If $p = 3$ then $G/K$ is one of the nearly--kaehler spaces
$$
Spin(8)/G_2 \text{ and } Spin(8)/SU(3).
$$
It is useful to note that
either the Dynkin diagram of $\gk$ is obtained by folding the diagram of 
$\gg$ as in \cite{dS2} --- all possibilities are listed in
the tables of \cite[pp. 245, 247]{dS2} --- or $\theta$ has form 
$\theta'\circ \Ad(g)$ 
where $\gg^{\theta'}$ is obtained by folding and $g \in G^{\theta'}$.
For example $Spin(8)/G_2$ is obtained by folding and $Spin(8)/SU(3)$
is derived from it as just described; and $E_6/F_4$ is obtained by folding
and $E_6/C_4$ is derived from it as just described.  For details of 
the latter see \cite[p. 291]{W1}.  We now run through that list.

\subsection{Case $G/K = SU(2m)/SO(2m)$.}  In order to find the 
restriction matrix
used by the LiE program, we consider the Cartan subalgebras

$\gs = \{ \diag\{u_1, \dots , u_m, -u_m, \dots ,-u_1\} \mid 
u_i \in \sqrt{-1}\,\R \}$ of $\gk$ and

$\gt = \{ \diag\{u_1, \dots , u_{2m}\}\mid u_i \in \sqrt{-1}\,\R, 
u_1 + \dots + u_{2m} = 0\}$ of $\gg$.

\noindent The simple roots of $\gg$ are the $\psi_i = \varepsilon_i -
\varepsilon_{i+1}$ for $1 \leqq i < 2m$, and the simple roots of $\gk$
are the $\varphi_i =  \varepsilon_i - \varepsilon_{i+1}$ 
(for $1 \leqq i < m$) and $\varphi_m = \varepsilon_{m-1} + \varepsilon_m$.
Thus the simple roots of $\gg$ have restriction to $\gs$ given by
$\psi_i|_\gs = \varphi_i$ and $\psi_{m+i}|_\gs = \varphi_{m-i}$ for
$1 \leqq i < m$ and $\psi_m|_\gs = 2\varepsilon_m = \varphi_m - \varphi_{m-1}$. 
Here is the relevant LiE routine {\tt branch\_A\_D.lie} for branching
from $SU(2m)$ to $SO(2m)$.  It takes arguments {\tt(m,v)}, where
$m > 0$ is an integer and {\tt v} is a vector of length $2m-1$ consisting
of non--negative integers, and branches
{\tt v}.  If {\tt m} is already set in LiE then only the argument 
{\tt v} is needed.
{\small
\begin{equation}\label{branch_A_D}
\begin{array}l
\verb! # file branch_A_D.lie                                            # ! \\
\verb! # usage: branch_A_D(m,v) branches v from SU(2m) to SO(2m)        # ! \\
\verb! # and branch_A_D(v) does the same if m is already defined in LiE # ! \\
\verb! branch_A_D(int m;vec v) = setdefault(Lie_group(1,2*m - 1)); ! \\
\verb! res_rt = null(2*m-1,m); ! \\
\verb! for i=1 to m do res_rt[i,i] = 1 od; ! \\
\verb! res_rt[m,m-1] = -1; res_rt[m,m] = 1; ! \\
\verb! for i=1 to m-1 do res_rt[m+i,m-i] = 1 od; ! \\
\verb! res_wt = i_Cartan*res_rt*Cartan(Lie_group(4,m))/det_Cartan; ! \\
\verb! answer = branch(v,Lie_group(4,m),res_wt); ! \\
\verb! print("the branching of "+v+" from SU("+2*m+") to SO("+2*m+") is"); ! \\
\verb! answer ! \\ 
\verb! branch_A_D(vec v) = setdefault(Lie_group(1,2*m - 1)); ! \\
\verb! res_rt = null(2*m-1,m); ! \\
\verb! for i=1 to m do res_rt[i,i] = 1 od; ! \\
\verb! res_rt[m,m-1] = -1; res_rt[m,m] = 1; ! \\
\verb! for i=1 to m-1 do res_rt[m+i,m-i] = 1 od; ! \\
\verb! res_wt = i_Cartan*res_rt*Cartan(Lie_group(4,m))/det_Cartan; ! \\
\verb! answer = branch(v,Lie_group(4,m),res_wt); ! \\
\verb! print("the branching of "+v+" from SU("+2*m+") to SO("+2*m+") is"); ! \\
\verb! answer !
\end{array}
\end{equation}}

Here is an example of its use:
{\small
\begin{quote}
\begin{verbatim}
read branch_A_D.lie 
branch_A_D(8,[1,1,0,0,0,0,0,0,0,0,0,0,0,0,0])
     the branching of [1,1,0,0,0,0,0,0,0,0,0,0,0,0,0] from SU(16) 
     to SO(16) is 1X[1,0,0,0,0,0,0,0] +1X[1,1,0,0,0,0,0,0]
m=5
branch_A_D([1,2,0,0,0,0,0,0,0])
     the branching of [1,2,0,0,0,0,0,0,0] from SU(10) to SO(10) is
     1X[1,0,0,0,0] +1X[1,1,0,0,0] +1X[1,2,0,0,0] +1X[3,0,0,0,0]
\end{verbatim}
\end{quote}}

\subsection{Case $G/K = SU(2m+1)/SO(2m+1)$.}  In order to find the 
restriction matrix we consider the Cartan subalgebras

$\gs = \{ \diag\{u_1, \dots , u_m, 0, -u_m, \dots ,-u_1\} \mid
u_i \in \sqrt{-1}\,\R \}$ of $\gk$ and

$\gt = \{ \diag\{u_1, \dots , u_{2m+1}\}\mid u_i \in \sqrt{-1}\,\R,
u_1 + \dots + u_{2m+1} = 0\}$ of $\gg$.

\noindent The simple roots of $\gg$ are the $\psi_i = \varepsilon_i -
\varepsilon_{i+1}$ for $1 \leqq i \leqq 2m$, and the simple roots of $\gk$
are the $\varphi_i =  \varepsilon_i - \varepsilon_{i+1}$
(for $1 \leqq i < m$) and $\varphi_m = \varepsilon_m$ (short simple root).
Thus the simple roots of $\gg$ have restriction to $\gs$ given by
$\psi_i|_\gs = \varphi_i$ and $\psi_{m+i}|_\gs = \varphi_{m+1-i}$ for
$1 \leqq i \leqq m$.  
Here is the relevant LiE routine {\tt branch\_A\_B.lie} for branching
from $SU(2m+1)$ to $SO(2m+1)$.  It takes arguments {\tt(m,v)}, where
$m > 0$ is an integer and {\tt v} is a vector of length $2m$ consisting
of non--negative integers, and branches
{\tt v}.  If {\tt m} is already set in LiE then only the argument
{\tt v} is needed.

{\small
\begin{equation}\label{branch_A_B}
\begin{array}l
\verb! # file branch_A_B.lie                                            # ! \\
\verb! # usage: branch_A_B(m,v) branches v from SU(2m+1) to SO(2m+1)    # ! \\
\verb! # and branch_A_B(v) does the same if m is already defined in LiE # ! \\
\verb! branch_A_B(int m;vec v) = setdefault(Lie_group(1,2*m)); ! \\
\verb! res_rt = null(2*m,m); ! \\
\verb! for i=1 to m do res_rt[i,i] = 1 od; ! \\
\verb! for i=1 to m do res_rt[m+i,m+1-i] = 1 od; ! \\
\verb! res_wt = i_Cartan*res_rt*Cartan(Lie_group(2,m))/det_Cartan; ! \\
\verb! answer = branch(v,Lie_group(2,m),res_wt); ! \\ 
\verb! print("the branching of "+v+" from SU("+(2*m+1)+") to \ ! \\
\verb!      SO("+(2*m+1)+") is"); ! \\
\verb! answer ! \\
\verb! branch_A_B(vec v) = setdefault(Lie_group(1,2*m)); ! \\
\verb! res_rt = null(2*m,m); ! \\
\verb! for i=1 to m do res_rt[i,i] = 1 od; ! \\
\verb! for i=1 to m do res_rt[m+i,m+1-i] = 1 od; ! \\
\verb! res_wt = i_Cartan*res_rt*Cartan(Lie_group(2,m))/det_Cartan; ! \\
\verb! answer = branch(v,Lie_group(2,m),res_wt); ! \\
\verb! print("the branching of "+v+" from SU("+(2*m+1)+") to \ ! \\
\verb!      SO("+(2*m+1)+") is"); ! \\
\verb! answer !
\end{array}
\end{equation}}

\subsection{Case $G/K = SU(2m)/Sp(m)$.}  This case is quite similar to the
case of $SU(2m)/SO(2m)$ above.  We consider the Cartan subalgebras

$\gs = \{ \diag\{u_1, \dots , u_m, -u_m, \dots ,-u_1\} \mid
u_i \in \sqrt{-1}\,\R \}$ of $\gk$ and

$\gt = \{ \diag\{u_1, \dots , u_{2m}\}\mid u_i \in \sqrt{-1}\,\R,
u_1 + \dots + u_{2m} = 0\}$ of $\gg$.

\noindent The simple roots of $\gg$ are the $\psi_i = \varepsilon_i -
\varepsilon_{i+1}$ for $1 \leqq i < 2m$, and the simple roots of $\gk$
are the $\varphi_i =  \varepsilon_i - \varepsilon_{i+1}$
(for $1 \leqq i < m$) and $\varphi_m = 2\varepsilon_m$.
Thus the simple roots of $\gg$ have restriction to $\gs$ given by
$\psi_i|_\gs = \varphi_i$ and $\psi_{m+i}|_\gs = \varphi_{m-i}$ for
$1 \leqq i < m$ and $\psi_m|_\gs = 2\varepsilon_m = \varphi_m$.
Here is the relevant LiE routine {\tt branch\_A\_C.lie} for branching
from $SU(2m)$ to $Sp(m)$.  It takes arguments {\tt(m,v)}, where
$m > 0$ is an integer and {\tt v} is a vector of length $2m-1$ consisting
of non--negative integers, and branches
{\tt v}.  If {\tt m} is already set in LiE then only the argument
{\tt v} is needed.

{\small
\begin{equation}\label{branch_A_C}
\begin{array}l
\verb!# file branch_A_C.lie                                            # ! \\
\verb!# usage: branch_A_C(m,v) branches v from SU(2m) to Sp(m)         # ! \\
\verb!# and branch_A_C(v) does the same if m is already defined in LiE # ! \\
\verb!branch_A_C(int m;vec v) = setdefault(Lie_group(1,2*m - 1)); ! \\
\verb!res_rt = null(2*m-1,m); ! \\
\verb!for i=1 to m-1 do res_rt[i,i] = 1 od; ! \\
\verb!res_rt[m,m] = 1; ! \\
\verb!for i=1 to m-1 do res_rt[m+i,m-i] = 1 od; ! \\
\verb!res_wt = i_Cartan*res_rt*Cartan(Lie_group(3,m))/det_Cartan; ! \\
\verb!answer = branch(v,Lie_group(3,m),res_wt); ! \\
\verb!print("the branching of "+v+" from SU("+2*m+") to Sp("+m+") is"); ! \\
\verb!answer ! \\ 
\verb!branch_A_C(vec v) = setdefault(Lie_group(1,2*m - 1)); ! \\
\verb!res_rt = null(2*m-1,m); ! \\
\verb!for i=1 to m-1 do res_rt[i,i] = 1 od; ! \\
\verb!res_rt[m,m] = 1; ! \\
\verb!for i=1 to m-1 do res_rt[m+i,m-i] = 1 od; ! \\
\verb!res_wt = i_Cartan*res_rt*Cartan(Lie_group(3,m))/det_Cartan; ! \\
\verb!answer = branch(v,Lie_group(3,m),res_wt); ! \\
\verb!print("the branching of "+v+" from SU("+2*m+") to Sp("+m+") is"); ! \\
\verb!answer !
\end{array}
\end{equation}}

\subsection{Cases $G/K = SO(2p + 2 + 2q)/\{SO(2p+1)\times SO(2q+1)\}$
($p,q$ not both $0$).}  We use the Cartan subalgebras

$\gt = \{u := \diag\{u_1,\dots u_{p+q+1}, - u_{p+q+1}, \dots -u_1\}\mid u_i
\in \sqrt{-1}\,\R\}$ and
$\gs = \{u \in \gt \mid u_{p+q+1} = 0\}$.

\noindent
Now $\gg$ has simple roots $\psi_i = \varepsilon_i - \varepsilon_{i+1}$ for
$1\leqq i \leqq p+q$ and $\psi_{p+q+1} = \varepsilon_{p+q} + 
\varepsilon_{p+q+1}$.  
The subalgebra $\gk$ has simple roots 
$\varphi_i = \varepsilon_i - \varepsilon_{i+1}$ for $1 \leqq i < p$ and for
$p+1 \leqq i < p+q$, $\varphi_p = \varepsilon_p$, and $\varphi_{p+q} = 
\varepsilon_{p+q}$.  
Thus the simple roots of $\gg$ have restriction to $\gs$ given by
$\psi_i|_\gs = \varphi_i$ for $1 \leqq i < p$ and $p < i < p+q$,
$\psi_p|_\gs = \varphi_p - \sum_1^q \varphi_{p+j}$,
$\psi_{p+q}|_\gs = \varphi_{p+q}$, and
$\psi_{p+q+1}|_\gs = \varphi_{p+q}$.
Here is the relevant LiE routine {\tt branch\_D\_BB.lie} for branching
from $SO(2p+2q+2)$ to $SO(2p+1) \times SO(2q+1)$.  It takes arguments 
{\tt(p,q,v)}, where $p,q \geqq 0$ are integers not both zero and {\tt v} 
is a vector of length $p+q+1$ consisting of non--negative integers, and 
branches {\tt v}.  If {\tt p} and {\tt q} are already set in LiE then only 
the argument {\tt v} is needed.

{\small
\begin{equation}\label{branch_D_BB}
\begin{array}l
\verb! # file branch_D_BB.lie                                     # ! \\
\verb! # usage: branch_D_BB(p,q,v) branches v from SO(2p+2q+2)    # ! \\
\verb! #        to SO(2p+1)xSO(2q+1), and branch_D_BB(v) does the # ! \\
\verb! #        same if p and q are already defined in LiE        # ! \\
\verb! branch_D_BB(int p,q;vec v) = setdefault(Lie_group(4,p+q+1)); ! \\
\verb! res_rt = null(p+q+1,p+q); ! \\
\verb! for i=1 to p-1 do res_rt[i,i] = 1 od; ! \\
\verb! res_rt[p,p] = 1; for j=p+1 to p+q do res_rt[p,j] = -1 od; ! \\
\verb! for i=p+1 to p+q-1 do res_rt[i,i] = 1 od; ! \\
\verb! res_rt[p+q,p+q] = 1; res_rt[p+q+1,p+q] = 1; ! \\
\verb! res_wt = i_Cartan*res_rt*Cartan(Lie_group(2,p)*Lie_group(2,q))/4; ! \\
\verb! answer = branch(v,Lie_group(2,p)*Lie_group(2,q),res_wt); ! \\
\verb! print("the branching of "+v+" from SO("+(2*p+2*q+2)+") to \ ! \\
\verb!      SO("+(2*p+1)+")xSO("+(2*q+1)+") is"); ! \\
\verb! answer ! \\ 
\verb! branch_D_BB(vec v) = setdefault(Lie_group(4,p+q+1)); ! \\
\verb! res_rt = null(p+q+1,p+q); ! \\
\verb! for i=1 to p-1 do res_rt[i,i] = 1 od; ! \\
\verb! res_rt[p,p] = 1; for j=p+1 to p+q do res_rt[p,j] = -1 od; ! \\
\verb! for i=p+1 to p+q-1 do res_rt[i,i] = 1 od; ! \\
\verb! res_rt[p+q,p+q] = 1; res_rt[p+q+1,p+q] = 1; ! \\
\verb! res_wt = i_Cartan*res_rt*Cartan(Lie_group(2,p)*Lie_group(2,q))/4; ! \\
\verb! answer = branch(v,Lie_group(2,p)*Lie_group(2,q),res_wt); ! \\
\verb! print("the branching of "+v+" from SO("+(2*p+2*q+2)+") to \ ! \\
\verb!      SO("+(2*p+1)+")xSO("+(2*q+1)+") is"); ! \\
\verb! answer !
\end{array}
\end{equation}}

\subsection{Case $G/K = Spin(8)/G_2$.}

In this case $\gg$ has simple root system $\{\psi_1,\psi_2,\psi_3,\psi_4\}$
numbered as in the Introduction, and $\gk$ has simple root system
$\{\varphi_1,\varphi_2\}$ as follows.  The Cartan subalgebra $\gs$ of $\gk$
is the subspace of the Cartan subalgebra $\gt$ of $\gg$ given by
$\psi_1 = \psi_3 = \psi_4$.  See \cite{dS1}.  The root restrictions
are $\psi_1|_\gs = \psi_3|_\gs = \psi_4|_\gs = \varphi_1$
(the short simple root of $\gk$) and $\psi_2|_\gs = \varphi_2$ (the long simple
root of $\gk$).  
Now the Lie routine for branching from $Spin(8)$ to $G_2$ is
{\small
\begin{equation}\label{branch_D4_G2}
\begin{array}l
\verb! # file branch_D4_G2.lie # ! \\
\verb! branch_D4_G2(vec v) = setdefault(D4); ! \\
\verb! res_rt = [[1,0],[0,1],[1,0],[1,0]]; ! \\
\verb! res_wt = i_Cartan*res_rt*Cartan(G2))/det_Cartan; ! \\
\verb! answer = branch(v,G2,res_wt); ! \\
\verb! print("the branching of "+v+" from Spin(8) to G2 is"); ! \\
\verb! answer ! 
\end{array}
\end{equation}}

\subsection{Cases $G/K = Spin(8)/SU(3)$, $G/K = E_6/F_4$ and $G/K = E_6/Sp(4)$.}
There the restriction matrices can be computed as in the case of 
$Spin(8)/G_2$, or one can use the small database of maximal subalgebras 
built into LiE.  That database is accessible by the commands

\centerline{{\tt res\_mat(A2,D4)}, {\tt res\_mat(F4,E6)} and {\tt res\_mat(C4,E6)}}

\noindent The corresponding LiE routines are
{\small
\begin{equation}\label{branch_D4_A2}
\begin{array}l
\verb!# file branch_D4_A2.lie #! \\
\verb!branch_D4_A2(vec v) = setdefault(D4);! \\
\verb!answer = branch(v,A2,res_mat(A2,D4));! \\
\verb!print("the branching of "+v+" from D4 to the triality A2 is"); ! \\
\verb!answer!
\end{array}
\end{equation}}
and
{\small
\begin{equation}\label{branch_E6_F4}
\begin{array}l
\verb!# file branch_E6_F4.lie #! \\
\verb!branch_E6_F4(vec v) = setdefault(E6);! \\
\verb!answer = branch(v,F4,res_mat(F4,E6));! \\
\verb!print("the branching of "+v+" from E6 to F4 is"); answer! 
\end{array}
\end{equation}}
and 
{\small
\begin{equation}\label{branch_E6_C4}
\begin{array}l
\verb!# file branch_E6_C4.lie #! \\
\verb!branch_E6_C4(vec v) = setdefault(E6);! \\
\verb!answer = branch(v,C4,res_mat(C4,E6));! \\
\verb!print("the branching of "+v+" from E6 to Sp(4) is"); answer! 
\end{array}
\end{equation}}

\section{Cases: $\gg$ is simple and $\gk$ is the centralizer of
a toral subalgebra}
\label{sec4}
\setcounter{equation}{0}

These cases were covered in Section \ref{sec2c}.  If $\gk$ is the
centralizer of a toral subalgebra of $\gg$ it is the fixed point set of
an automorphism.  For if $G$ is a connected
Lie group with Lie algebra $\gg$ and $K$ is the analytic subgroup for
$\gk$ we can choose $g \in G$ such that the powers 
$\{g^n \mid n \in \Z\}$ form a dense subgroup of the identity component
of the center of $K$, and then $\gk$ is the fixed point set of $\Ad(g)$.

\section{Cases: $\gg$ is simple, $\theta^2 = 1$ and $\rank \gk = \rank \gg$}
\label{sec5}
\setcounter{equation}{0}

In this section we apply the method of Section \ref{sec2d} and 
run through the cases where $n_r = 2$, i.e.\ the cases
where $G/K$ is a riemannian symmetric space.  The other cases of equal
rank will be considered in the next section.

It turns out that, in the classical group symmetric space cases, we do 
not have to renumber the roots of $\Psi_\gs$,
while the renumbering is needed for most of the exceptional cases.

\subsection{Case  $(G,K) = (SO(2p +2q+1),SO(2p)\times SO(2q+1))$
where $p \geqq 2$ and $q \geqq 0$.}\label{B_DB1}
Here $w_\gs$ reverses
the order of $\psi_1, \dots , \psi_{r-1}$ but does not move $\psi_i$
for $i > r$, so $w_\gs(\beta_\gg)$ attaches to the diagram of $\gs$
at $\psi_{r-2}$.  Now  the simple roots of $\gk$ in Bourbaki root order are
given by
$\{\psi_1, \dots , \psi_{r-1}, w_\gs(\beta_\gg); \psi_{r+1}, \dots , \psi_n\}$. 
Here is the LiE program {\tt branch\_B\_DB.lie} that branches the representation
of $SO(2p+2q+1)$, specified by a vector {\tt v} of length $p+q$, to 
$SO(2p) \times SO(2q+1)$.  Usage is {\tt branch\_B\_DB(p,q,v)}, but 
if $p$ and $q$ are already set in LiE one can just use {\tt branch\_B\_DB(v)}.
{\small
\begin{equation}\label{branch_B_DB}
\begin{array}l
\verb!#file branch_B_DB.lie                                                     #!\\
\verb!#usage: branch_B_DB(p,q,v) branches v from SO(2p+2q+1) to SO(2p)xSO(2q+1) #!\\
\verb!#and branch_B_DB(v) does the same if p and q are already defined in LiE   #!\\
\verb! branch_B_DB(int p,q; vec v) = setdefault(Lie_group(2,p+q)); ! \\
\verb! loc JJ = A1*A1; loc KK = A1; loc LL = Lie_group(4,p+q); ! \\
\verb! u = null(p+q-1); for i = 1 to p+q-1 do u[i]=i od; ! \\
\verb! if q > 0 then for i = p to p+q-1 do u[i]=i+1 od fi; ! \\
\verb! ws = reduce(long_word^r_reduce(long_word,u)); ! \\
\verb! RR = id(p+q); RR[p] = W_rt_action(high_root,ws); ! \\
\verb! if p >= 3 then JJ = Lie_group(4,p) fi; ! \\
\verb! if q >= 2 then KK = Lie_group(2,q) fi; ! \\
\verb! if q > 0 then LL = JJ*KK fi; answer = branch(v,LL,res_mat(RR)); ! \\
\verb! print("the branching of "+v+" from SO("+(2*p+2*q+1)+") " ! \\
\verb!    "to SO("+2*p+")xSO("+(2*q+1)+") is "); ! \\
\verb! answer ! \\
\verb! branch_B_DB(vec v) = setdefault(Lie_group(2,p+q)); ! \\
\verb! loc JJ = A1*A1; loc KK = A1; loc LL = Lie_group(4,p+q); ! \\
\verb! u = null(p+q-1); for i = 1 to p+q-1 do u[i]=i od; ! \\
\verb! if q > 0 then for i = p to p+q-1 do u[i]=i+1 od fi; ! \\
\verb! ws = reduce(long_word^r_reduce(long_word,u)); ! \\
\verb! RR = id(p+q); RR[p] = W_rt_action(high_root,ws); ! \\
\verb! if p >= 3 then JJ = Lie_group(4,p) fi; ! \\
\verb! if q >= 2 then KK = Lie_group(2,q) fi; ! \\
\verb! if q > 0 then LL = JJ*KK fi; answer = branch(v,LL,res_mat(RR)); ! \\
\verb! print("the branching of "+v+" from SO("+(2*p+2*q+1)+") " ! \\
\verb!    "to SO("+2*p+")xSO("+(2*q+1)+") is "); ! \\
\verb! answer ! 
\end{array}
\end{equation}
}

\subsection{Case $(G,K) = (SO(2p+2q),SO(2p)\times SO(2q))$
where $p, q \geqq 2$.}\label{D_DD}
Again $w_\gs$ reverses
the order of $\psi_1, \dots , \psi_{r-1}$ but does not move $\psi_i$
for $r < i \leqq n-2$, so $w_\gs(\beta_\gg)$ attaches to the diagram of $\gs$
at $\psi_{r-2}$ (or doesn't attach, if $r = 2$).  Now  the simple roots 
of $\gk$ in Bourbaki root order are
$\{\psi_1, \dots , \psi_{r-1}, w_\gs(\beta_\gg); \psi_{r+1}, \dots , \psi_n\}$.
Here is the LiE program for restriction of representations from
$SO(2p+2q)$ to $SO(2p)\times SO(2q)$.
{\small
\begin{equation}\label{branch_D_DD}
\begin{array}l
\verb! #file branch_D_DD.lie                                                   # ! \\
\verb! #usage: branch_D_DD(p,q,v) branches v from SO(2p+2q) to SO(2p)xSO(2q)   # ! \\
\verb! #and branch_D_DD(v) does the same if p and q are already defined in LiE # ! \\
\verb! branch_D_DD(int p,q; vec v) = setdefault(Lie_group(4,p+q)); ! \\
\verb! loc JJ = A1*A1; loc KK = A1*A1; ! \\
\verb! u = null(p+q-1); for i = 1 to p+q-1 do u[i]=i od; ! \\
\verb! for i = p to p+q-1 do u[i]=i+1 od; ! \\
\verb! ws = reduce(long_word^r_reduce(long_word,u)); ! \\
\verb! RR = id(p+q); RR[p] = W_rt_action(high_root,ws); ! \\
\verb! if p >= 3 then JJ = Lie_group(4,p) fi; ! \\
\verb! if q >= 3 then KK = Lie_group(4,q) fi; ! \\
\verb! answer = branch(v,JJ*KK,res_mat(RR)); ! \\
\verb! print("the branching of "+v+" from SO("+(2*p+2*q)+") " ! \\
\verb! "to SO("+2*p+")xSO("+(2*q)+") is "); ! \\
\verb! answer ! \\
\verb! branch_D_DD(vec v) = setdefault(Lie_group(4,p+q)); ! \\
\verb! loc JJ = A1*A1; loc KK = A1*A1; ! \\
\verb! u = null(p+q-1); for i = 1 to p+q-1 do u[i]=i od; ! \\
\verb! for i = p to p+q-1 do u[i]=i+1 od; ! \\
\verb! ws = reduce(long_word^r_reduce(long_word,u)); ! \\
\verb! RR = id(p+q); RR[p] = W_rt_action(high_root,ws); ! \\
\verb! if p >= 3 then JJ = Lie_group(4,p) fi; ! \\
\verb! if q >= 3 then KK = Lie_group(4,q) fi; ! \\
\verb! answer = branch(v,JJ*KK,res_mat(RR)); ! \\
\verb! print("the branching of "+v+" from SO("+(2*p+2*q)+") " ! \\
\verb! "to SO("+2*p+")xSO("+(2*q)+") is "); ! \\
\verb! answer ! 
\end{array}
\end{equation}
}

\subsection{Case $(G,K) = (Sp(p+q),Sp(p)\times Sp(q))$
where $p, q \geqq 1$.}\label{C_CC}
Again $w_\gs$ reverses
the order of $\psi_1, \dots , \psi_{r-1}$ but does not move $\psi_i$
for $i > r$, so $w_\gs(\beta_\gg)$ attaches to the diagram of $\gs$
at $\psi_{r-1}$ (or doesn't attach, if $r = 1$).  Now  the simple roots
of $\gk$ in Bourbaki root order are
$\{\psi_1, \dots , \psi_{r-1}, w_\gs(\beta_\gg); \psi_{r+1}, \dots , \psi_n\}$.
The LiE program for restriction of representations from $Sp(p+q)$ to
$Sp(p)\times Sp(q)$ is 
{\small
\begin{equation}\label{branch_C_CC}
\begin{array}l
\verb! #file branch_C_CC.lie                                              # ! \\
\verb! #usage: branch_C_CC(p,q,v) branches v from Sp(p+q) to Sp(p)xSp(q), # ! \\
\verb! #branch_C_CC(v) does the same if p, q are already defined in LiE   # ! \\
\verb! branch_C_CC(int p,q; vec v) = setdefault(Lie_group(3,p+q)); ! \\
\verb! loc JJ = A1; loc KK = A1; ! \\
\verb! u = null(p+q-1); if p == 1 then for i = 1 to q do u[i]=i+1 od fi; ! \\
\verb! if p >= 2 then for i = 1 to p+q-1 do u[i]=i od fi; ! \\
\verb! if p >= 2 then for i = p to p+q-1 do u[i]=i+1 od fi; ! \\
\verb! ws = reduce(long_word^r_reduce(long_word,u)); ! \\
\verb! RR = id(p+q); RR[p] = W_rt_action(high_root,ws); ! \\
\verb! if p >= 2 then JJ = Lie_group(3,p) fi; ! \\
\verb! if q >= 2 then KK = Lie_group(3,q) fi; ! \\
\verb! answer = branch(v,JJ*KK,res_mat(RR)); ! \\
\verb! print("the branching of "+v+" from Sp("+(p+q)+") " ! \\
\verb! "to Sp("+p+")xSp("+q+") is "); ! \\
\verb! answer ! \\
\verb! branch_C_CC(vec v) = setdefault(Lie_group(3,p+q)); ! \\
\verb! loc JJ = A1; loc KK = A1; ! \\
\verb! u = null(p+q-1); if p == 1 then for i = 1 to q do u[i]=i+1 od fi; ! \\
\verb! if p >= 2 then for i = 1 to p+q-1 do u[i]=i od fi; ! \\
\verb! if p >= 2 then for i = p to p+q-1 do u[i]=i+1 od fi; ! \\
\verb! ws = reduce(long_word^r_reduce(long_word,u)); ! \\
\verb! RR = id(p+q); RR[p] = W_rt_action(high_root,ws); ! \\
\verb! if p >= 2 then JJ = Lie_group(3,p) fi; ! \\
\verb! if q >= 2 then KK = Lie_group(3,q) fi; ! \\
\verb! answer = branch(v,JJ*KK,res_mat(RR)); ! \\
\verb! print("the branching of "+v+" from Sp("+(p+q)+") " ! \\
\verb! "to Sp("+p+")xSp("+q+") is "); ! \\
\verb! answer ! 
\end{array}
\end{equation}}

\subsection{Case $G_2/A_1A_1$.} \label{G2}
Now we run through the exceptional group cases.  First suppose that
$G = G_2$.  Then $K = A_1A_1$ with simple roots
\setlength{\unitlength}{.6 mm}
\begin{picture}(20,10)
\put(4,1){\circle{3}}
\put(5.5,1){\line(1,0){6.6}}
\put(15,1){\makebox(0,0){\circle{3}}}
\put(5.0,-.3){\line(1,0){8.0}}
\put(5.0,2.3){\line(1,0){8.0}}
\put(9,1){\makebox(0,0){$\langle$}}
\put(5,4){\makebox(0,0)[b]{$\psi_1$}}
\put(15,4){\makebox(0,0)[b]{$\psi_2$}}
\end{picture}.
Here $\beta_\gg = 3\psi_1 + 2\psi_2$ and $\gamma = \psi_2$,
$\Psi_\gs = \{\psi_1\}$, and $w_\gs(\beta_\gg) = \beta_\gg$.
Thus the LiE program (if one wants to bother with it in this case) is
{\small
\begin{equation}\label{branch_G2_A1A1}
\begin{array}l
\verb!# file branch_G2_A1A1.lie # ! \\
\verb!branch_G2_A1A1(vec v) = setdefault(G2); ! \\
\verb!ws = reduce(long_word^r_reduce(long_word,[1])); ! \\
\verb!RR = id(2); RR[2] = W_rt_action(high_root,ws); ! \\
\verb!answer = branch(v,A1A1,res_mat(RR)); ! \\
\verb!print("the branching of "+v+" from G2 to A1A1 is"); answer !
\end{array}
\end{equation}}

\subsection{Cases $F_4/A_1C_3$ and $F_4/B_4$.}\label{F4}
Next suppose that $G = F_4$.  The simple root system is
\setlength{\unitlength}{.6 mm}
\begin{picture}(45,10)
\put(5,1){\circle{2}}
\put(6,1){\line(1,0){8}}
\put(15,1){\circle{2}}
\put(15.2,0.2){\line(1,0){9.6}}
\put(15.2,1.8){\line(1,0){9.6}}
\put(25,1){\circle{2}}
\put(26,1){\line(1,0){9}}
\put(37,1){\makebox(0,0){\circle{2}}}
\put(20,1){\makebox(0,0){$\rangle$}}
\put(5,4){\makebox(0,0)[b]{$\psi_1$}}
\put(15,4){\makebox(0,0)[b]{$\psi_2$}}
\put(25,4){\makebox(0,0)[b]{$\psi_3$}}
\put(35,4){\makebox(0,0)[b]{$\psi_4$}}
\end{picture}
and the negative of the 
maximal root $\beta_\gg = 2\psi_1 + 3\psi_2 + 4\psi_3 + 2\psi_4$
attaches at $\psi_1$. Thus there are two possibilities:
$\gamma = \psi_1$ and $K = A_1C_3$, or $\gamma = \psi_4$ and $K = B_4$.
The corresponding LiE programs are given by
{\small
\begin{equation}\label{branch_F4_A1C3}
\begin{array}l
\verb!# file branch_F4_A1C3.lie # ! \\
\verb!branch_F4_A1C3(vec v) = setdefault(F4); ! \\
\verb!ws = reduce(long_word^r_reduce(long_word,[2,3,4])); ! \\
\verb!RR = id(4); RR[1] = W_rt_action(high_root,ws); ! \\
\verb!RR[2] = id(4)[4]; RR[4] = id(4)[2]; ! \\
\verb!answer = branch(v,A1C3,res_mat(RR)); ! \\
\verb!print("the branching of "+v+" from F4 to A1C3 is"); answer !
\end{array}
\end{equation}}
and
{\small
\begin{equation}\label{branch_F4_B4}
\begin{array}l
\verb!# file branch_F4_B4.lie # ! \\
\verb!branch_F4_B4(vec v) = setdefault(F4); ! \\
\verb!ws = reduce(long_word^r_reduce(long_word,[1,2,3])); ! \\
\verb!RR = id(4); RR[1] = W_rt_action(high_root,ws); ! \\
\verb!RR[2] = id(4)[1]; RR[3] = id(4)[2]; RR[4] = id(4)[3]; ! \\
\verb!answer = branch(v,B4,res_mat(RR)); ! \\
\verb!print("the branching of "+v+" from F4 to B4 is"); answer !
\end{array}
\end{equation}}

\subsection{Cases $E_6/A_1A_5$ and $E_6/A_5A_1$.}\label{E6}
Now suppose that $G = E_6$.  Then there are three simple roots of coefficient 
$2$ in the maximal root.  All of them differ by automorphisms of the extended
Dynkin diagram, so the corresponding subalgebras $\gk$ differ by an
automorphism, but in \cite{EW} we will need to distinguish between them ,
so we treat them
separately.  The simple root system is
\setlength{\unitlength}{.6 mm}
\begin{picture}(55,16)(0,-5)
\put(5,3){\circle{2}}
\put(6,3){\line(1,0){8}}
\put(15,3){\circle{2}}
\put(16,3){\line(1,0){8}}
\put(25,3){\circle{2}}
\put(26,3){\line(1,0){8}}
\put(35,3){\circle{2}}
\put(36,3){\line(1,0){8}}
\put(45,3){\circle{2}}
\put(25,2){\line(0,-1){8}}
\put(25,-7){\circle{2}}
\put(5,6){\makebox(0,0)[b]{$\psi_1$}}
\put(15,6){\makebox(0,0)[b]{$\psi_3$}}
\put(25,6){\makebox(0,0)[b]{$\psi_4$}}
\put(35,6){\makebox(0,0)[b]{$\psi_5$}}
\put(45,6){\makebox(0,0)[b]{$\psi_6$}}
\put(29,-6){\makebox(0,0)[l]{$\psi_2$}}
\end{picture}
and the negative of the 
maximal root $\beta_\gg = \psi_1 + 2\psi_2 + 2\psi_3 + 3\psi_4
+ 2\psi_5 + \psi_6$ attaches at $\psi_2$.  There are three equivalent
possibilities:
(i) $\gamma = \psi_3$ and $K = A_1A_5$, (ii) $\gamma = \psi_5$ and 
$K = A_5A_1$, and (iii) $\gamma = \psi_2$ and $K = A_5A_1$.
For the first, the LiE routine is
{\small
\begin{equation}\label{branch_E6_A1A5}
\begin{array}l
\verb!# file branch_E6_A1A5.lie # ! \\
\verb!branch_E6_A1A5(vec v) = setdefault(E6); ! \\
\verb!ws = reduce(long_word^r_reduce(long_word,[1,2,4,5,6])); ! \\
\verb!RR = id(6); RR[6] = W_rt_action(high_root,ws); ! \\
\verb!RR[3] = id(6)[4]; RR[4] = id(6)[5]; RR[5] = id(6)[6]; ! \\
\verb!answer = branch(v,A1A5,res_mat(RR)); ! \\
\verb!print("the branching of "+v+" from E6 to A1A5 is"); answer !
\end{array}
\end{equation}}
The second is quite similar,
{\small
\begin{equation}\label{branch_E6_A5A1}
\begin{array}l
\verb!# file branch_E6_A5A1.lie # ! \\
\verb!branch_E6_A5A1(vec v) = setdefault(E6); ! \\
\verb!ws = reduce(long_word^r_reduce(long_word,[1,2,3,4,6])); ! \\
\verb!RR = id(6); RR[1] = W_rt_action(high_root,ws); ! \\
\verb!RR[2] = id(6)[1]; RR[5] = id(6)[2]; ! \\
\verb!answer = branch(v,A5A1,res_mat(RR)); ! \\
\verb!print("the branching of "+v+" from E6 to A5A1 is"); answer !
\end{array}
\end{equation}}
and the third is a bit different,
{\small
\begin{equation}\label{branch_E6_A5A1a}
\begin{array}l
\verb!# file branch_E6_A5A1a.lie # ! \\
\verb!branch_E6_A5A1a(vec v) = setdefault(E6); ! \\
\verb!ws = reduce(long_word^r_reduce(long_word,[1,3,4,5,6])); ! \\
\verb!RR = id(6); RR[6] = W_rt_action(high_root,ws); ! \\
\verb!RR[2] = id(6)[3]; RR[3] = id(6)[4]; RR[4] = id(6)[5]; RR[5] = id(6)[6]; ! \\
\verb!answer = branch(v,A5A1,res_mat(RR)); ! \\
\verb!print("the branching of "+v+" from E6 to A5A1 is"); answer !
\end{array}
\end{equation}}

\subsection{Cases $E_7/A_1D_6$, $E_7/D_6A_1$ and $E_7/A_7$.}\label{E7}
Next, let $G = E_7$.  Then there are three simple roots of coefficient $2$
in the maximal root. Two of them differ by an automorphism of the extended
Dynkin diagram, so the corresponding subalgebras $\gk$ differ by an
automorphism, but we will need to distinguish between them in \cite{EW}.
So, as in some of the $E_6$ cases above, we treat them
separately.  The simple root system is
\setlength{\unitlength}{.6 mm}
\begin{picture}(62,16)(0,-5)
\put(5,3){\circle{2}}
\put(6,3){\line(1,0){8}}
\put(15,3){\circle{2}}
\put(16,3){\line(1,0){8}}
\put(25,3){\circle{2}}
\put(26,3){\line(1,0){8}}
\put(35,3){\circle{2}}
\put(36,3){\line(1,0){8}}
\put(45,3){\circle{2}}
\put(46,3){\line(1,0){9}}
\put(56,3){\circle{2}}
\put(25,2){\line(0,-1){8}}
\put(25,-7){\circle{2}}
\put(5,6){\makebox(0,0)[b]{$\psi_1$}}
\put(15,6){\makebox(0,0)[b]{$\psi_3$}}
\put(25,6){\makebox(0,0)[b]{$\psi_4$}}
\put(35,6){\makebox(0,0)[b]{$\psi_5$}}
\put(45,6){\makebox(0,0)[b]{$\psi_6$}}
\put(55,6){\makebox(0,0)[b]{$\psi_7$}}
\put(27,-4){\makebox(0,0)[l]{$\psi_2$}}
\end{picture}
and the negative of the 
maximal root $\beta_\gg = 2\psi_1 + 2\psi_2 + 3\psi_3 + 4\psi_4
+ 3\psi_5 + 2\psi_6 + \psi_7$ attaches at $\psi_1$.  The
possibilities are (i) $\gamma = \psi_1$ and $K = A_1D_6$, 
(ii) $\gamma = \psi_6$ and $K = D_6A_1$,
and (iii) $\gamma = \psi_2$ and $K = A_7$.  For first of these the LiE
routine is
{\small
\begin{equation}\label{branch_E7_A1D6}
\begin{array}l
\verb!# file branch_E7_A1D6.lie # ! \\
\verb!branch_E7_A1D6(vec v) = setdefault(E7); ! \\
\verb!ws = reduce(long_word^r_reduce(long_word,[2,3,4,5,6,7])); ! \\
\verb!RR = id(7); RR[1] = W_rt_action(high_root,ws); ! \\
\verb!RR[2] = id(7)[7]; RR[3] = id(7)[6]; RR[4] = id(7)[5]; ! \\
\verb!RR[5] = id(7)[4]; RR[6] = id(7)[3]; RR[7] = id(7)[2]; ! \\
\verb!answer = branch(v,A1D6,res_mat(RR)); ! \\
\verb!print("the branching of "+v+" from E7 to A1D6 is"); answer !
\end{array}
\end{equation}}
For obvious reasons the second is similar
{\small
\begin{equation}\label{branch_E7_D6A1}
\begin{array}l
\verb!# file branch_E7_D6A1.lie # ! \\
\verb!branch_E7_D6A1(vec v) = setdefault(E7); ! \\
\verb!ws = reduce(long_word^r_reduce(long_word,[1,2,3,4,5,7])); ! \\
\verb!RR = id(7); RR[1] =  W_rt_action(high_root,ws); ! \\
\verb!RR[2] = id(7)[1]; RR[6] = id(7)[2]; ! \\
\verb!answer = branch(v,D6A1,res_mat(RR)); ! \\
\verb!print("the branching of "+v+" from E7 to D6A1 is"); answer !
\end{array}
\end{equation}}
and the third is a bit different
{\small
\begin{equation}\label{branch_E7_A7}
\begin{array}l
\verb!# file branch_E7_A7.lie # ! \\
\verb!branch_E7_A7(vec v) = setdefault(E7); ! \\
\verb!ws = reduce(long_word^r_reduce(long_word,[1,3,4,5,6,7])); ! \\
\verb!RR = id(7); RR[7] = W_rt_action(high_root,ws); ! \\
\verb!RR[2] = id(7)[3]; RR[3] = id(7)[4]; RR[4] = id(7)[5]; ! \\
\verb!RR[5] = id(7)[6]; RR[6] = id(7)[7]; ! \\
\verb!answer = branch(v,A7,res_mat(RR)); ! \\
\verb!print("the branching of "+v+" from E7 to A7 is"); answer !
\end{array}
\end{equation}}

\subsection{Cases $E_8/D_8$ and $E_8/E_7A_1$}\label{E8}
Finally, suppose that $G = E_8$.  The simple root system is
\setlength{\unitlength}{.6 mm}
\begin{picture}(67,10)(0,-5)
\put(5,3){\circle{2}}
\put(6,3){\line(1,0){8}}
\put(15,3){\circle{2}}
\put(16,3){\line(1,0){8}}
\put(25,3){\circle{2}}
\put(26,3){\line(1,0){8}}
\put(35,3){\circle{2}}
\put(36,3){\line(1,0){8}}
\put(45,3){\circle{2}}
\put(46,3){\line(1,0){8}}
\put(55,3){\circle{2}}
\put(56,3){\line(1,0){9}}
\put(66,3){\circle{2}}
\put(25,2){\line(0,-1){8}}
\put(25,-7){\circle{2}}
\put(5,6){\makebox(0,0)[b]{$\psi_1$}}
\put(15,6){\makebox(0,0)[b]{$\psi_3$}}
\put(25,6){\makebox(0,0)[b]{$\psi_4$}}
\put(35,6){\makebox(0,0)[b]{$\psi_5$}}
\put(45,6){\makebox(0,0)[b]{$\psi_6$}}
\put(55,6){\makebox(0,0)[b]{$\psi_7$}}
\put(65,6){\makebox(0,0)[b]{$\psi_8$}}
\put(27,-4){\makebox(0,0)[l]{$\psi_2$}}
\end{picture}
and the negative of the 
maximal root $\beta_\gg = 2\psi_1 + 3\psi_2 + 4\psi_3 + 6\psi_4 +
5\psi_5 + 4\psi_6 + 3\psi_7 + 2\psi_8$ attaches at $\psi_8$.  The 
possibilities are (i) $\gamma = \psi_1$ and $K = D_8$ and 
(ii) $\gamma = \psi_8$ and $K = E_7A_1$.  In the first case the LiE routine
is
{\small
\begin{equation}\label{branch_E8_D8}
\begin{array}l
\verb!# file branch_E8_D8.lie # ! \\
\verb!branch_E8_D8(vec v) = setdefault(E8); ! \\
\verb!ws = reduce(long_word^r_reduce(long_word,[2,3,4,5,6,7,8])); ! \\
\verb!RR = id(8); RR[1] = W_rt_action(high_root,ws); ! \\
\verb!RR[2] = id(8)[8]; RR[3] = id(8)[7]; RR[4] = id(8)[6]; ! \\
\verb!RR[6] = id(8)[4]; RR[7] = id(8)[3]; RR[8] = id(8)[2]; ! \\
\verb!answer = branch(v,D8,res_mat(RR)); ! \\
\verb!print("the branching of "+v+" from E8 to D8 is"); answer !
\end{array}
\end{equation}}
and in the second it is
{\small
\begin{equation}\label{branch_E8_E7A1}
\begin{array}l
\verb!# file branch_E8_E7A1.lie # ! \\
\verb!branch_E8_E7A1(vec v) = setdefault(E8); ! \\
\verb!ws = reduce(long_word^r_reduce(long_word,[1,2,3,4,5,6,7])); ! \\
\verb!RR = id(8); RR[8] = W_rt_action(high_root,ws); ! \\
\verb!answer = branch(v,E7A1,res_mat(RR)); ! \\
\verb!print("the branching of "+v+" from E8 to E7A1 is"); answer !
\end{array}
\end{equation}}

Taking into account the results of Sections \ref{sec3} and \ref{sec4},
now for every compact connected, simply connected symmetric space $G/K$, 
$G$ simple, we have shown how to compute the restriction to $K$ of any 
irreducible finite dimensional representation of $G$.  The irreducible
compact connected, simply connected symmetric spaces $G/K$ with $\gg$
not simple are the are the simply connected simple Lie group manifolds,
which were covered in Section \ref{sec2a}

\section{Cases: $\gg$ is simple, $\theta^3 = 1$ or $\theta^5 = 1$, 
and $\rank \gk = \rank \gg$}
\label{sec6}
\setcounter{equation}{0}

The maximal connected subgroups of maximal rank in a compact connected
Lie group were described by A. Borel and J. de Siebenthal \cite{BdS}.
Most of them are symmetric subgroups, and their classification can be used
in the classification of symmetric spaces \cite{W1}.  The ones that
are symmetric were considered in Section \ref{sec5}.  The others correspond
to the simple roots $\gamma$ whose coefficient in the maximal root $\beta_\gg$
is an odd prime, necessarily $3$ or $5$.  The ones for prime $3$
are given by $(G,K) = (G_2,A_2), (F_4,A_2A_2), 
(E_6,A_2A_2A_2), (E_7,A_2A_5), (E_7,A_5A_2),
(E_8,A_8), \text{ and } (E_8,E_6A_2)$.  The one for prime $5$ is given by
$(G,K) = (E_8,A_4A_4)$.  In all cases a simple root system for $\gk$ is
given by $\Psi_\gk = (\Psi_\gg \setminus \{\gamma\})\cup \{-\beta_\gg\}$,
so we can use the methods of Section \ref{sec2c} as in Section \ref{sec5}.

In all of these cases one can rely on a LiE database to produce a
restriction matrix {\tt res\_mat(Y,X)}.  However the applications in
\cite{EW} require that we keep track of which root of {\tt Y} comes from
which root of {\tt X}, and LiE scrambles the root order, so we generally 
have to do this by hand.  In each case we indicate which elements of
the simple root system $\Phi = \{\varphi_1, \dots , \varphi_n\}$ of $\gk$
come from which elements of 
$(\Psi \setminus \{\gamma\})\cup \{w_\gs(\beta_\gg)\}$, where we
try to use the least complicated correspondence.

\subsection{Case $G/K = G_2/A_2$.}

Here $\gamma = \psi_1$, $\gk$ is of type $A_2 = \gs\gu(3)$, and 
$\Phi = \{\varphi_1,\varphi_2\}$ where 
$\varphi_1$ comes from $w_\gs(\beta_\gg)$ and
$\varphi_2$ comes from $\psi_2$.
This is indicated in terms of the Dynkin diagrams, by

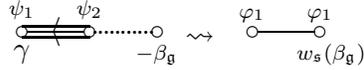
\centerline{
\setlength{\unitlength}{.75 mm}
\begin{picture}(34,10)(0,-2)
\put(5,0){\circle{2}}
\put(5,4){\makebox(0,0){$\scriptstyle\psi_1$}}
\put(5,-4){\makebox(0,0){$\gamma$}}
\put(5,-.95){\line(1,0){12}}
\put(6,0){\line(1,0){10}}
\put(5,.95){\line(1,0){12}}
\put(11,0){\makebox(0,0){$\langle$}}
\put(17,0){\circle{2}}
\put(17,4){\makebox(0,0){$\scriptstyle\psi_2$}}
\multiput(18,0)(1,0){11}{\makebox(0,0){.}}
\put(29,0){\circle{2}}
\put(29,-4){\makebox(0,0){$\scriptstyle-\beta_\gg$}}
\end{picture}$\rightsquigarrow$
\begin{picture}(22,10)(0,-2)
\put(5,0){\circle{2}}
\put(6,0){\line(1,0){10}}
\put(17,0){\circle{2}}
\put(5,3){\makebox(0,0){$\scriptstyle\varphi_1$}}
\put(17,3){\makebox(0,0){$\scriptstyle\varphi_1$}}
\put(19,-4){\makebox(0,0){$\scriptstyle w_\gs(\beta_\gg)$}}
\end{picture}
}
\noindent The corresponding LiE routine is
{\small
\begin{equation}\label{branch_G2_A2}
\begin{array}l
\verb! # file branch_G2_A2.lie # ! \\
\verb! branch_G2_A2(vec v) = setdefault(G2); !\\
\verb! ws = reduce(long_word^r_reduce(long_word,[2])); !\\
\verb! RR = id(2); RR[1] = W_rt_action(high_root,ws); !\\
\verb! answer = branch(v,A2,res_mat(RR)); !\\
\verb! print("the branching of "+v+" from G2 to A2 is "); answer !
\end{array}
\end{equation}}

\subsection{Case $G/K = F_4/A_2A_2$.}

Here $\gamma = \psi_2$, $\gk$ is of type $A_2A_2$, and $\Phi = 
\{\varphi_1, \varphi_2, \varphi_3, \varphi_4\}$ where 
$\varphi_1$ comes from $\psi_1$,
$\varphi_2$ comes from $w_\gs(\beta_\gg)$,
$\varphi_3$ comes from $\psi_3$, and
$\varphi_4$ comes from $\psi_4$.  This is indicated in terms of 
Dynkin diagrams by

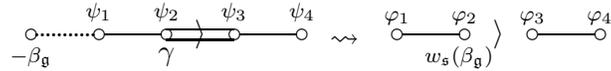
\centerline{
\setlength{\unitlength}{.75 mm}
\begin{picture}(58,10)(0,-2)    
\put(5,0){\circle{2}}
\put(17,0){\circle{2}}
\put(29,0){\circle{2}}
\put(41,0){\circle{2}}
\put(53,0){\circle{2}}
\multiput(6,0)(1,0){11}{\makebox(0,0){.}}
\put(18,0){\line(1,0){10}}
\put(29,-.95){\line(1,0){12}}
\put(29,.95){\line(1,0){12}}
\put(35,0){\makebox(0,0){$\rangle$}}
\put(42,0){\line(1,0){10}}
\put(17,4){\makebox(0,0){$\scriptstyle\psi_1$}}
\put(29,4){\makebox(0,0){$\scriptstyle\psi_2$}}
\put(41,4){\makebox(0,0){$\scriptstyle\psi_3$}}
\put(53,4){\makebox(0,0){$\scriptstyle\psi_4$}}
\put(29,-4){\makebox(0,0){$\gamma$}}
\put(5,-4){\makebox(0,0){$\scriptstyle-\beta_\gg$}}
\end{picture}$\rightsquigarrow$
\begin{picture}(46,10)(0,-2)    
\put(5,0){\circle{2}}
\put(17,0){\circle{2}}
\put(29,0){\circle{2}}
\put(41,0){\circle{2}}
\put(6,0){\line(1,0){10}}
\put(23,0){\makebox(0,0){$\rangle$}}
\put(30,0){\line(1,0){10}}
\put(5,3){\makebox(0,0){$\scriptstyle\varphi_1$}}
\put(17,3){\makebox(0,0){$\scriptstyle\varphi_2$}}
\put(29,3){\makebox(0,0){$\scriptstyle\varphi_3$}}
\put(41,3){\makebox(0,0){$\scriptstyle\varphi_4$}}
\put(16,-4){\makebox(0,0){$\scriptstyle w_\gs(\beta_\gg)$}}
\end{picture}}

\noindent The corresponding LiE routine is
{\small
\begin{equation}\label{branch_F4_A2A2}
\begin{array}l
\verb! # file branch_F4_A2A2.lie # ! \\
\verb! branch_F4_A2A2(vec v) = setdefault(F4); ! \\
\verb! ws = reduce(long_word^r_reduce(long_word,[1,3,4])); ! \\
\verb! RR = id(4); RR[2] = W_rt_action(high_root,ws); ! \\
\verb! answer = branch(v,A2A2,res_mat(RR)); ! \\
\verb! print("the branching of "+v+" from F4 to A2A2 is "); answer ! 
\end{array}
\end{equation}}

\subsection{Case $G/K = E_6/A_2A_2A_2$.}

Here $\gamma = \psi_4$, $\gk$ is of type $A_2A_2A_2$, and $\Phi =
\{\varphi_1, \varphi_2, \varphi_3, \varphi_4, \varphi_5, \varphi_6\}$ where
$\varphi_1$ comes from $\psi_1$,
$\varphi_2$ comes from $\psi_3$,
$\varphi_3$ comes from $\psi_2$,
$\varphi_4$ comes from $w_\gs(\beta_\gg)$,
$\varphi_5$ comes from $\psi_5$, and
$\varphi_6$ comes from $\psi_6$.
This is indicated in terms of Dynkin diagrams by

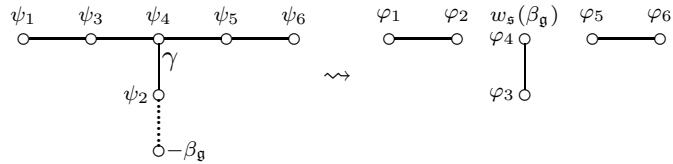
\centerline{
\setlength{\unitlength}{.75 mm}
\begin{picture}(58,30)(0,-22)    
\put(5,0){\circle{2}}
\put(17,0){\circle{2}}
\put(29,0){\circle{2}}
\put(41,0){\circle{2}}
\put(53,0){\circle{2}}
\put(29,-10){\circle{2}}
\put(29,-20){\circle{2}}
\put(6,0){\line(1,0){10}}
\put(18,0){\line(1,0){10}}
\put(30,0){\line(1,0){10}}
\put(42,0){\line(1,0){10}}
\put(29,-1){\line(0,-1){8}}
\multiput(29,-11)(0,-1){9}{\makebox(0,0){.}}
\put(5,4){\makebox(0,0){$\scriptstyle\psi_1$}}
\put(17,4){\makebox(0,0){$\scriptstyle\psi_3$}}
\put(29,4){\makebox(0,0){$\scriptstyle\psi_4$}}
\put(41,4){\makebox(0,0){$\scriptstyle\psi_5$}}
\put(53,4){\makebox(0,0){$\scriptstyle\psi_6$}}
\put(25,-10){\makebox(0,0){$\scriptstyle\psi_2$}}
\put(31,-4){\makebox(0,0){$\gamma$}}
\put(34,-20){\makebox(0,0){$\scriptstyle-\beta_\gg$}}
\end{picture}\raisebox{30pt}{$\rightsquigarrow$}
\begin{picture}(58,30)(0,-22)    
\put(5,0){\circle{2}}
\put(17,0){\circle{2}}
\put(29,0){\circle{2}}
\put(41,0){\circle{2}}
\put(53,0){\circle{2}}
\put(29,-10){\circle{2}}
\put(6,0){\line(1,0){10}}
\put(42,0){\line(1,0){10}}
\put(29,-1){\line(0,-1){8}}
\put(5,4){\makebox(0,0){$\scriptstyle\varphi_1$}}
\put(17,4){\makebox(0,0){$\scriptstyle\varphi_2$}}
\put(25,0){\makebox(0,0){$\scriptstyle\varphi_4$}}
\put(41,4){\makebox(0,0){$\scriptstyle\varphi_5$}}
\put(53,4){\makebox(0,0){$\scriptstyle\varphi_6$}}
\put(25,-10){\makebox(0,0){$\scriptstyle\varphi_3$}}
\put(29,4){\makebox(0,0){$\scriptstyle w_\gs(\beta_\gg)$}}
\end{picture}
}

\noindent The corresponding LiE routine is
{\small
\begin{equation}\label{branch_E6_A2A2A2}
\begin{array}l
\verb! # file branch_E6_A2A2A2.lie # ! \\
\verb! branch_E6_A2A2A2(vec v) = setdefault(E6); ! \\
\verb! ws = reduce(long_word^r_reduce(long_word,[1,2,3,5,6])); ! \\
\verb! RR = id(6); RR[4] = W_rt_action(high_root,ws); ! \\
\verb! RR[2] = id(6)[3]; RR[3] = id(6)[2];  ! \\
\verb! answer = branch(v,A2A2A2,res_mat(RR)); ! \\ 
\verb! print("the branching of "+v+" from E6 to A2A2A2 is"); answer !
\end{array}
\end{equation}}

\subsection{Case $G/K = E_7/A_2A_5$.}

Here $\gamma = \psi_3$, $\gk$ is of type $A_2A_5$, and $\Phi =
\{\varphi_1,\varphi_2,\varphi_3,\varphi_4,\varphi_5,\varphi_6,\varphi_7\}$
where
$\varphi_1$ comes from $\psi_1$,
$\varphi_2$ comes from $w_\gs(\beta_\gg)$,
$\varphi_3$ comes from $\psi_2$,
$\varphi_4$ comes from $\psi_4$,
$\varphi_5$ comes from $\psi_5$,
$\varphi_6$ comes from $\psi_6$,
and
$\varphi_7$ comes from $\psi_7$.
This is indicated in terms of Dynkin diagrams by

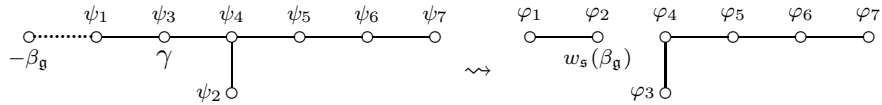
\centerline{
\setlength{\unitlength}{.75 mm}
\begin{picture}(82,20)(0,-12)    
\put(5,0){\circle{2}}
\put(17,0){\circle{2}}
\put(29,0){\circle{2}}
\put(41,0){\circle{2}}
\put(53,0){\circle{2}}
\put(65,0){\circle{2}}
\put(77,0){\circle{2}}
\put(41,-10){\circle{2}}
\put(18,0){\line(1,0){10}}
\put(30,0){\line(1,0){10}}
\put(42,0){\line(1,0){10}}
\put(54,0){\line(1,0){10}}
\put(66,0){\line(1,0){10}}
\put(41,-1){\line(0,-1){8}}
\multiput(6,0)(1,0){11}{\makebox(0,0){.}}
\put(5,-4){\makebox(0,0){$\scriptstyle-\beta_\gg$}}
\put(17,4){\makebox(0,0){$\scriptstyle\psi_1$}}
\put(29,4){\makebox(0,0){$\scriptstyle\psi_3$}}
\put(41,4){\makebox(0,0){$\scriptstyle\psi_4$}}
\put(53,4){\makebox(0,0){$\scriptstyle\psi_5$}}
\put(65,4){\makebox(0,0){$\scriptstyle\psi_6$}}
\put(77,4){\makebox(0,0){$\scriptstyle\psi_7$}}
\put(37,-10){\makebox(0,0){$\scriptstyle\psi_2$}}
\put(29,-4){\makebox(0,0){$\gamma$}}
\end{picture}\raisebox{10pt}{$\rightsquigarrow$}
\begin{picture}(70,20)(12,-12)    
\put(17,0){\circle{2}}
\put(29,0){\circle{2}}
\put(41,0){\circle{2}}
\put(53,0){\circle{2}}
\put(65,0){\circle{2}}
\put(77,0){\circle{2}}
\put(41,-10){\circle{2}}
\put(18,0){\line(1,0){10}}
\put(42,0){\line(1,0){10}}
\put(54,0){\line(1,0){10}}
\put(66,0){\line(1,0){10}}
\put(41,-1){\line(0,-1){8}}
\put(17,4){\makebox(0,0){$\scriptstyle\varphi_1$}}
\put(29,4){\makebox(0,0){$\scriptstyle\varphi_2$}}
\put(41,4){\makebox(0,0){$\scriptstyle\varphi_4$}}
\put(53,4){\makebox(0,0){$\scriptstyle\varphi_5$}}
\put(65,4){\makebox(0,0){$\scriptstyle\varphi_6$}}
\put(77,4){\makebox(0,0){$\scriptstyle\varphi_7$}}
\put(37,-10){\makebox(0,0){$\scriptstyle\varphi_3$}}
\put(29,-4){\makebox(0,0){$\scriptstyle w_\gs(\beta_\gg)$}}
\end{picture}
}

\noindent The corresponding LiE routine is
{\small
\begin{equation}\label{branch_E7_A2A5}
\begin{array}l
\verb! # file branch_E7_A2A5.lie # ! \\ 
\verb! branch_E7_A2A5(vec v) = setdefault(E7); ! \\
\verb! ws = reduce(long_word^r_reduce(long_word,[1,2,4,5,6,7])); ! \\
\verb! RR = id(7); RR[2] = W_rt_action(high_root,ws); RR[3] = id(7)[2]; ! \\
\verb! res = res_mat(RR); ! \\
\verb! answer = branch(v,A2A5,res_mat(RR)); ! \\
\verb! print("the branching of "+v+" from E7 to A2A5 is"); answer ! 
\end{array}
\end{equation}}

\subsection{Case $G/K = E_7/A_5A_2$.}

Here $\gamma = \psi_5$, $\gk$ is of type $A_5A_2$, and $\Phi =
\{\varphi_1,\varphi_2,\varphi_3,\varphi_4,\varphi_5,\varphi_6,\varphi_7\}$
where
$\varphi_1$ comes from $\psi_1$,
$\varphi_2$ comes from $\psi_3$,
$\varphi_3$ comes from $\psi_4$,
$\varphi_4$ comes from $\psi_2$,
$\varphi_5$ comes from $w_\gs(\beta_\gg)$,
$\varphi_6$ comes from $\psi_6$,
and
$\varphi_7$ comes from $\psi_7$.
This is indicated in terms of Dynkin diagrams by

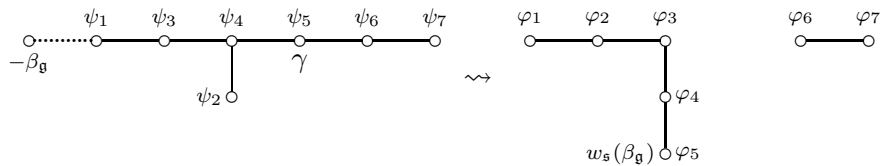
\centerline{
\setlength{\unitlength}{.75 mm}
\begin{picture}(82,30)(0,-22)    
\put(5,0){\circle{2}}
\put(17,0){\circle{2}}
\put(29,0){\circle{2}}
\put(41,0){\circle{2}}
\put(53,0){\circle{2}}
\put(65,0){\circle{2}}
\put(77,0){\circle{2}}
\put(41,-10){\circle{2}}
\put(18,0){\line(1,0){10}}
\put(30,0){\line(1,0){10}}
\put(42,0){\line(1,0){10}}
\put(54,0){\line(1,0){10}}
\put(66,0){\line(1,0){10}}
\put(41,-1){\line(0,-1){8}}
\multiput(6,0)(1,0){11}{\makebox(0,0){.}}
\put(5,-4){\makebox(0,0){$\scriptstyle-\beta_\gg$}}
\put(17,4){\makebox(0,0){$\scriptstyle\psi_1$}}
\put(29,4){\makebox(0,0){$\scriptstyle\psi_3$}}
\put(41,4){\makebox(0,0){$\scriptstyle\psi_4$}}
\put(53,4){\makebox(0,0){$\scriptstyle\psi_5$}}
\put(65,4){\makebox(0,0){$\scriptstyle\psi_6$}}
\put(77,4){\makebox(0,0){$\scriptstyle\psi_7$}}
\put(37,-10){\makebox(0,0){$\scriptstyle\psi_2$}}
\put(53,-4){\makebox(0,0){$\gamma$}}
\end{picture}\raisebox{30pt}{$\rightsquigarrow$}
\begin{picture}(70,30)(12,-22)    
\put(17,0){\circle{2}}
\put(29,0){\circle{2}}
\put(41,0){\circle{2}}
\put(65,0){\circle{2}}
\put(77,0){\circle{2}}
\put(41,-10){\circle{2}}
\put(41,-20){\circle{2}}
\put(18,0){\line(1,0){10}}
\put(30,0){\line(1,0){10}}
\put(66,0){\line(1,0){10}}
\put(41,-1){\line(0,-1){8}}
\put(41,-11){\line(0,-1){8}}
\put(17,4){\makebox(0,0){$\scriptstyle\varphi_1$}}
\put(29,4){\makebox(0,0){$\scriptstyle\varphi_2$}}
\put(41,4){\makebox(0,0){$\scriptstyle\varphi_3$}}
\put(65,4){\makebox(0,0){$\scriptstyle\varphi_6$}}
\put(77,4){\makebox(0,0){$\scriptstyle\varphi_7$}}
\put(45,-10){\makebox(0,0){$\scriptstyle\varphi_4$}}
\put(45,-20){\makebox(0,0){$\scriptstyle\varphi_5$}}
\put(33,-20){\makebox(0,0){$\scriptstyle w_\gs(\beta_\gg)$}}
\end{picture}
}

\noindent The corresponding LiE routine is
{\small
\begin{equation}\label{branch_E7_A5A2}
\begin{array}l
\verb! # file branch_E7_A5A2.lie # ! \\
\verb! branch_E7_A5A2(vec v) = setdefault(E7); ! \\
\verb! ws = reduce(long_word^r_reduce(long_word,[1,2,3,4,6,7])); ! \\
\verb! RR = id(7); RR[5] = W_rt_action(high_root,ws);  ! \\
\verb! RR[2] = id(7)[3]; RR[3] = id(7)[4]; RR[4] = id(7)[2]; ! \\
\verb! answer = branch(v,A5A2,res_mat(RR)); ! \\
\verb! print("the branching of "+v+" from E7 to A5A2 is"); answer ! 
\end{array}
\end{equation}}

\subsection{Case $G/K = E_8/A_8$.} 

Here $\gamma = \psi_2$, $\gk$ is of type $A_8$, and $\Phi =
\{\varphi_1,\varphi_2,\varphi_3,\varphi_4,\varphi_5,\varphi_6,\varphi_7,
\varphi_8\}$ where
$\varphi_1$ comes from $w_\gs(\beta_\gg)$,
$\varphi_2$ comes from $\psi_1$,
$\varphi_3$ comes from $\psi_3$,
$\varphi_4$ comes from $\psi_4$,
$\varphi_5$ comes from $\psi_5$,
$\varphi_6$ comes from $\psi_6$,
$\varphi_7$ comes from $\psi_7$,
and
$\varphi_8$ comes from $\psi_8$.
This is indicated in terms of Dynkin diagrams by

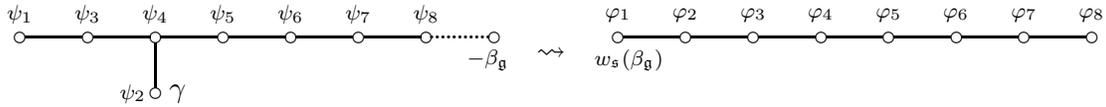
\centerline{
\setlength{\unitlength}{.75 mm}
\begin{picture}(94,20)(0,-12)    
\put(5,0){\circle{2}}
\put(17,0){\circle{2}}
\put(29,0){\circle{2}}
\put(41,0){\circle{2}}
\put(53,0){\circle{2}}
\put(65,0){\circle{2}}
\put(77,0){\circle{2}}
\put(89,0){\circle{2}}
\put(29,-10){\circle{2}}
\put(6,0){\line(1,0){10}}
\put(18,0){\line(1,0){10}}
\put(30,0){\line(1,0){10}}
\put(42,0){\line(1,0){10}}
\put(54,0){\line(1,0){10}}
\put(66,0){\line(1,0){10}}
\put(29,-1){\line(0,-1){8}}
\multiput(78,0)(1,0){11}{\makebox(0,0){.}}
\put(88,-4){\makebox(0,0){$\scriptstyle-\beta_\gg$}}
\put(5,4){\makebox(0,0){$\scriptstyle\psi_1$}}
\put(17,4){\makebox(0,0){$\scriptstyle\psi_3$}}
\put(29,4){\makebox(0,0){$\scriptstyle\psi_4$}}
\put(41,4){\makebox(0,0){$\scriptstyle\psi_5$}}
\put(53,4){\makebox(0,0){$\scriptstyle\psi_6$}}
\put(65,4){\makebox(0,0){$\scriptstyle\psi_7$}}
\put(77,4){\makebox(0,0){$\scriptstyle\psi_8$}}
\put(25,-10){\makebox(0,0){$\scriptstyle\psi_2$}}
\put(33,-10){\makebox(0,0){$\gamma$}}
\end{picture}\raisebox{17pt}{\enskip$\rightsquigarrow$\enskip}
\begin{picture}(94,20)(0,-12)    
\put(5,0){\circle{2}}
\put(17,0){\circle{2}}
\put(29,0){\circle{2}}
\put(41,0){\circle{2}}
\put(53,0){\circle{2}}
\put(65,0){\circle{2}}
\put(77,0){\circle{2}}
\put(89,0){\circle{2}}
\put(6,0){\line(1,0){10}}
\put(18,0){\line(1,0){10}}
\put(30,0){\line(1,0){10}}
\put(42,0){\line(1,0){10}}
\put(54,0){\line(1,0){10}}
\put(66,0){\line(1,0){10}}
\put(78,0){\line(1,0){10}}
\put(7,-4){\makebox(0,0){$\scriptstyle w_\gs(\beta_\gg)$}}
\put(5,4){\makebox(0,0){$\scriptstyle\varphi_1$}}
\put(17,4){\makebox(0,0){$\scriptstyle\varphi_2$}}
\put(29,4){\makebox(0,0){$\scriptstyle\varphi_3$}}
\put(41,4){\makebox(0,0){$\scriptstyle\varphi_4$}}
\put(53,4){\makebox(0,0){$\scriptstyle\varphi_5$}}
\put(65,4){\makebox(0,0){$\scriptstyle\varphi_6$}}
\put(77,4){\makebox(0,0){$\scriptstyle\varphi_7$}}
\put(89,4){\makebox(0,0){$\scriptstyle\varphi_8$}}
\end{picture}
}

\noindent The corresponding LiE routine is
{\small
\begin{equation}\label{branch_E8_A8}
\begin{array}l
\verb! # file branch_E8_A8.lie # ! \\
\verb! branch_E8_A8(vec v) = setdefault(E8); ! \\
\verb! ws = reduce(long_word^r_reduce(long_word,[1,3,4,5,6,7,8])); ! \\
\verb! RR = id(8); RR[1] = W_rt_action(high_root,ws); RR[2] = id(8)[1];  ! \\
\verb! answer = branch(v,A8,res_mat(RR)); ! \\
\verb! print("the branching of "+v+" from E8 to A8 is"); answer !
\end{array}
\end{equation}}

\subsection{Case $G/K = E_8/E_6A_2$.}

Here $\gamma = \psi_7$, $\gk$ is of type $E_6A_2$, and $\Phi =
\{\varphi_1,\varphi_2,\varphi_3,\varphi_4,\varphi_5,\varphi_6,\varphi_7,
\varphi_8\}$ where
$\varphi_1$ comes from $\psi_1$,
$\varphi_2$ comes from $\psi_2$,
$\varphi_3$ comes from $\psi_3$,
$\varphi_4$ comes from $\psi_4$,
$\varphi_5$ comes from $\psi_5$,
$\varphi_6$ comes from $\psi_6$,
$\varphi_7$ comes from $w_\gs(\beta_\gg)$,
and
$\varphi_8$ comes from $\psi_8$.
This is indicated in terms of Dynkin diagrams by

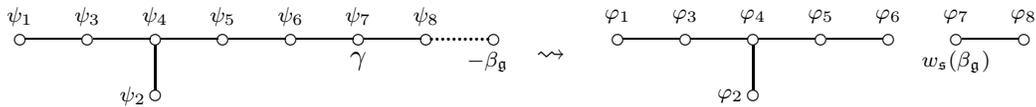
\centerline{
\setlength{\unitlength}{.75 mm}
\begin{picture}(94,20)(0,-12)    
\put(5,0){\circle{2}}
\put(17,0){\circle{2}}
\put(29,0){\circle{2}}
\put(41,0){\circle{2}}
\put(53,0){\circle{2}}
\put(65,0){\circle{2}}
\put(77,0){\circle{2}}
\put(89,0){\circle{2}}
\put(29,-10){\circle{2}}
\put(6,0){\line(1,0){10}}
\put(18,0){\line(1,0){10}}
\put(30,0){\line(1,0){10}}
\put(42,0){\line(1,0){10}}
\put(54,0){\line(1,0){10}}
\put(66,0){\line(1,0){10}}
\put(29,-1){\line(0,-1){8}}
\multiput(78,0)(1,0){11}{\makebox(0,0){.}}
\put(88,-4){\makebox(0,0){$\scriptstyle-\beta_\gg$}}
\put(5,4){\makebox(0,0){$\scriptstyle\psi_1$}}
\put(17,4){\makebox(0,0){$\scriptstyle\psi_3$}}
\put(29,4){\makebox(0,0){$\scriptstyle\psi_4$}}
\put(41,4){\makebox(0,0){$\scriptstyle\psi_5$}}
\put(53,4){\makebox(0,0){$\scriptstyle\psi_6$}}
\put(65,4){\makebox(0,0){$\scriptstyle\psi_7$}}
\put(77,4){\makebox(0,0){$\scriptstyle\psi_8$}}
\put(25,-10){\makebox(0,0){$\scriptstyle\psi_2$}}
\put(65,-4){\makebox(0,0){$\gamma$}}
\end{picture}\raisebox{17pt}{\enskip$\rightsquigarrow$\enskip}
\begin{picture}(82,20)(0,-12)    
\put(5,0){\circle{2}}
\put(17,0){\circle{2}}
\put(29,0){\circle{2}}
\put(41,0){\circle{2}}
\put(53,0){\circle{2}}
\put(65,0){\circle{2}}
\put(77,0){\circle{2}}
\put(29,-10){\circle{2}}
\put(6,0){\line(1,0){10}}
\put(18,0){\line(1,0){10}}
\put(30,0){\line(1,0){10}}
\put(42,0){\line(1,0){10}}
\put(66,0){\line(1,0){10}}
\put(29,-1){\line(0,-1){8}}
\put(65,-4){\makebox(0,0){$\scriptstyle w_\gs(\beta_\gg)$}}
\put(5,4){\makebox(0,0){$\scriptstyle\varphi_1$}}
\put(17,4){\makebox(0,0){$\scriptstyle\varphi_3$}}
\put(29,4){\makebox(0,0){$\scriptstyle\varphi_4$}}
\put(41,4){\makebox(0,0){$\scriptstyle\varphi_5$}}
\put(53,4){\makebox(0,0){$\scriptstyle\varphi_6$}}
\put(65,4){\makebox(0,0){$\scriptstyle\varphi_7$}}
\put(77,4){\makebox(0,0){$\scriptstyle\varphi_8$}}
\put(25,-10){\makebox(0,0){$\scriptstyle\varphi_2$}}
\end{picture}
}

\noindent The corresponding LiE routine is
{\small
\begin{equation}\label{branch_E8_E6A2}
\begin{array}l
\verb! # file branch_E8_E6A2.lie # ! \\
\verb! branch_E8_E6A2(vec v) = setdefault(E8); ! \\
\verb! ws = reduce(long_word^r_reduce(long_word,[1,2,3,4,5,6,8])); ! \\
\verb! RR = id(8); RR[7] = W_rt_action(high_root,ws); ! \\
\verb! answer = branch(v,E6A2,res_mat(RR)); ! \\
\verb! print("the branching of "+v+" from E8 to E6A2 is"); answer ! 
\end{array}
\end{equation}}

\subsection{Case $G/K = E_8/A_4A_4$.}

Here $\gamma = \psi_5$, $\gk$ is of type $A_4A_4$, and $\Phi =
\{\varphi_1,\varphi_2,\varphi_3,\varphi_4,\varphi_5,\varphi_6,\varphi_7,
\varphi_8\}$ where
$\varphi_1$ comes from $\psi_1$,
$\varphi_2$ comes from $\psi_3$,
$\varphi_3$ comes from $\psi_4$,
$\varphi_4$ comes from $\psi_2$,
$\varphi_5$ comes from $w_\gs(\beta_\gg)$,
$\varphi_6$ comes from $\psi_6$,
$\varphi_7$ comes from $\psi_7$,
and
$\varphi_8$ comes from $\psi_8$.
This is indicated in terms of Dynkin diagrams by

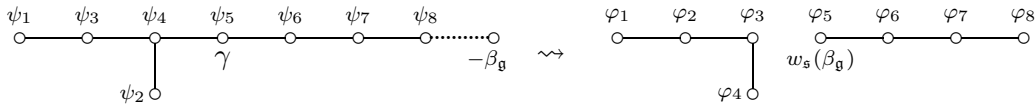
\centerline{
\setlength{\unitlength}{.75 mm}
\begin{picture}(94,20)(0,-12)    
\put(5,0){\circle{2}}
\put(17,0){\circle{2}}
\put(29,0){\circle{2}}
\put(41,0){\circle{2}}
\put(53,0){\circle{2}}
\put(65,0){\circle{2}}
\put(77,0){\circle{2}}
\put(89,0){\circle{2}}
\put(29,-10){\circle{2}}
\put(6,0){\line(1,0){10}}
\put(18,0){\line(1,0){10}}
\put(30,0){\line(1,0){10}}
\put(42,0){\line(1,0){10}}
\put(54,0){\line(1,0){10}}
\put(66,0){\line(1,0){10}}
\put(29,-1){\line(0,-1){8}}
\multiput(78,0)(1,0){11}{\makebox(0,0){.}}
\put(88,-4){\makebox(0,0){$\scriptstyle-\beta_\gg$}}
\put(5,4){\makebox(0,0){$\scriptstyle\psi_1$}}
\put(17,4){\makebox(0,0){$\scriptstyle\psi_3$}}
\put(29,4){\makebox(0,0){$\scriptstyle\psi_4$}}
\put(41,4){\makebox(0,0){$\scriptstyle\psi_5$}}
\put(53,4){\makebox(0,0){$\scriptstyle\psi_6$}}
\put(65,4){\makebox(0,0){$\scriptstyle\psi_7$}}
\put(77,4){\makebox(0,0){$\scriptstyle\psi_8$}}
\put(25,-10){\makebox(0,0){$\scriptstyle\psi_2$}}
\put(41,-4){\makebox(0,0){$\gamma$}}
\end{picture}\raisebox{17pt}{\enskip$\rightsquigarrow$\enskip}
\begin{picture}(82,20)(0,-12)    
\put(5,0){\circle{2}}
\put(17,0){\circle{2}}
\put(29,0){\circle{2}}
\put(41,0){\circle{2}}
\put(53,0){\circle{2}}
\put(65,0){\circle{2}}
\put(77,0){\circle{2}}
\put(29,-10){\circle{2}}
\put(6,0){\line(1,0){10}}
\put(18,0){\line(1,0){10}}
\put(42,0){\line(1,0){10}}
\put(54,0){\line(1,0){10}}
\put(66,0){\line(1,0){10}}
\put(29,-1){\line(0,-1){8}}
\put(41,-4){\makebox(0,0){$\scriptstyle w_\gs(\beta_\gg)$}}
\put(5,4){\makebox(0,0){$\scriptstyle\varphi_1$}}
\put(17,4){\makebox(0,0){$\scriptstyle\varphi_2$}}
\put(29,4){\makebox(0,0){$\scriptstyle\varphi_3$}}
\put(41,4){\makebox(0,0){$\scriptstyle\varphi_5$}}
\put(53,4){\makebox(0,0){$\scriptstyle\varphi_6$}}
\put(65,4){\makebox(0,0){$\scriptstyle\varphi_7$}}
\put(77,4){\makebox(0,0){$\scriptstyle\varphi_8$}}
\put(25,-10){\makebox(0,0){$\scriptstyle\varphi_4$}}
\end{picture}
}

\noindent The corresponding LiE routine is
{\small
\begin{equation}\label{branch_E8_A4A4}
\begin{array}l
\verb! # file branch_E8_A4A4.lie # ! \\
\verb! branch_E8_A4A4(vec v) = setdefault(E8); ! \\
\verb! ws = reduce(long_word^r_reduce(long_word,[1,2,3,4,6,7,8])); ! \\
\verb! RR = id(8); RR[5] = W_rt_action(high_root,ws); ! \\
\verb! RR[2] = id(8)[3]; RR[3] = id(8)[4]; RR[4] = id(8)[2]; ! \\
\verb! answer = branch(v,A4A4,res_mat(RR)); ! \\
\verb! print("the branching of "+v+" from E8 to A4A4 is"); answer ! 
\end{array}
\end{equation}}

This completes our branching project as described in Section 2.

\vskip 1 cm

\centerline{\begin{tabular}{ll}
MGE: & JAW: \\
Mathematical Sciences Institute & Department of Mathematics \\
Australian National University  & University of California \\
ACT 0200, Australia & Berkeley, California
94720--3840, U.S.A. \\
                                &                                      \\
{\tt meastwoo@member.ams.org} & {\tt jawolf@math.berkeley.edu}
\end{tabular}}

\end{document}